\documentclass{agujournal}

\usepackage{graphicx}
\usepackage{amsmath}
\usepackage{amssymb}
\usepackage{bm}
\usepackage{enumerate}   
\usepackage{algorithm}
\usepackage[noend]{algpseudocode}

\usepackage{xcolor}
\usepackage[colorinlistoftodos]{todonotes}
\usepackage[normalem]{ulem}

\journalname{Water Resources Research}

\begin{document}

\title{Randomized Truncated SVD Levenberg-Marquardt Approach to Geothermal Natural State and History Matching}

\authors{Elvar K. Bjarkason\affil{1},
		Oliver J. Maclaren\affil{1},
		John P. O'Sullivan\affil{1}, and
        Michael J. O'Sullivan\affil{1}}

\affiliation{1}{Department of Engineering Science, The University of Auckland, Auckland, New Zealand}

\correspondingauthor{Elvar Bjarkason}{ebja558@aucklanduni.ac.nz}
\correspondingauthor{Michael O'Sullivan}{m.osullivan@auckland.ac.nz}

%\begin{keypoints}
%\item Applied randomized TSVD methods to speed up parameter updating during inversion
%\item Randomization enables parallel solution of all adjoint and direct problems used for model updates
%\item The randomized approaches are shown to be considerably more efficient than standard approaches
%\end{keypoints}

\begin{abstract}
The Levenberg-Marquardt (LM) method is commonly used for inverting models used to describe geothermal, groundwater, or oil and gas reservoirs. In previous studies LM parameter updates have been made tractable for highly parameterized inverse problems with large data sets by applying matrix factorization methods or iterative linear solvers to approximately solve the update equations. 

Some studies have shown that basing model updates on the truncated singular value decomposition (TSVD) of a dimensionless sensitivity matrix achieved using Lanczos iteration can speed up the inversion of reservoir models. Lanczos iterations only require the sensitivity matrix times a vector and its transpose times a vector, which are found efficiently using adjoint and direct simulations without the expense of forming a large sensitivity matrix.

Nevertheless, Lanczos iteration has the drawback of being a serial process, requiring a separate adjoint solve and direct solve every Lanczos iteration. Randomized methods, developed for low-rank matrix approximation of large matrices, are more efficient alternatives to the standard Lanczos method. Here we develop LM variants which use randomized methods to find a TSVD of a dimensionless sensitivity matrix when updating parameters. The randomized approach offers improved efficiency by enabling simultaneous solution of all adjoint and direct problems for a parameter update.

\end{abstract}

\section{Introduction}
Inversion of a reservoir model is often performed by solving the following generalized least-squares problem: given $N_d$ observations ${\bm d}_{\text{obs}}$, minimize a regularized least-squares objective function, depending on the $N_m$ model parameters $\bm{m}$, which can be written as 
\begin{linenomath*}
\begin{equation} \label{eq:objective}
	\Phi ({\bm m}) = \Phi_{\text{d}} ({\bm m}) + \mu \Phi_{\text{m}} ({\bm m})   \, .
\end{equation} \end{linenomath*}
Here 
\begin{linenomath*}
\begin{equation} \label{eq:obsmismatch}
	\Phi_{\text{d}} ({\bm m}) = \left[ {\bm d} ({\bm m}) - {\bm d}_{\text{obs}}  \right]^T {\bm \Gamma}^{-1}_{\text{d}}   \left[ {\bm d} ({\bm m}) - {\bm d}_{\text{obs}}  \right]
\end{equation} \end{linenomath*}
is the observation mismatch term, where ${\bm d} ({\bm m})$ are simulated observations and ${\bm \Gamma}_{\text{d}}$ is the covariance matrix of the measurement noise. The last term in Eq. \eqref{eq:objective} is a regularizing or penalty function multiplied by a positive regularization weight $\mu$. Here we make the common assumption that the model penalty term can be written as
\begin{linenomath*}
\begin{equation} \label{eq:modelmismatch}
	\Phi_{\text{m}} ({\bm m}) = \left[ {\bm m} - {\bm m}_{\text{pr}}  \right]^T {\bm R}   \left[ {\bm m} - {\bm m}_{\text{pr}}  \right] \, ,
\end{equation} \end{linenomath*}
where ${\bm m}_{\text{pr}}$ is the prior or initial guess for the model parameters and ${\bm R}$ is a positive definite matrix, which imposes the correlation structure and preferred parameter conditions for the model parameters. When the model parameters are adequately described as being multivariate Gaussian then ${\bm R}$ can, for instance, be specified as the inverse of the prior covariance matrix \citep{oliver2008,oliver2011review}.

There are many methods available for minimizing \eqref{eq:objective} \citep{brochu2010,conn2009,nocedal2006,oliver2011review}. Here we only consider the Levenberg-Marquardt (LM) approach.

The LM method is the most commonly used method within the geothermal community for automatic inversion of geothermal reservoir models. It can be readily applied to geothermal problems using popular tools such as iTOUGH2 \citep{finsterle2007} and PEST \citep{doherty2016}. The LM approach was popularized for geothermal modeling with the advent of iTOUGH2 \citep{finsterle1995toughworkshop,finsterle1995solving,finsterle1997application}. PEST had been widely used for hydrological modeling \citep{zhou2014inverse,healy2010}, before being applied to geothermal problems \citep{austria2015,colina2013,osullivan2016orakeikorako}.

A drawback of standard LM implementations, such as those found in iTOUGH2 or PEST, is the need to evaluate ${\bm S}$, the $N_d$ by $N_m$ sensitivity matrix. In these two codes ${\bm S}$ is evaluated using expensive finite difference approximations, requiring $N_m + 1$ forward simulations. The evaluation of ${\bm S}$ can be made faster using direct or adjoint methods \citep{anterion1989,bjarkason2016,li2003,hinze2009,oliver2008,rodrigues2006adjoint}, but even this approach becomes infeasible for a large number of parameters and observations. 

This paper looks at accelerating inversion of reservoir parameters by applying a modified LM method. This new approach advances the work of \citet{tavakoli2010}, by applying recently developed methods from randomized linear algebra.

\citet{tavakoli2010} proposed basing LM model updates on an approximate truncated singular value decomposition (TSVD) of a dimensionless sensitivity matrix ${\bm S}_\text{D}$ found using Lanczos iteration. Their approach reduces computational cost for large $N_d$ and $N_m$ by avoiding forming the large matrix ${\bm S}$. Instead their method uses adjoint and direct simulations to evaluate ${\bm S}$ and ${\bm S}^T$ times vectors. 
However, since the Lanczos approach is iterative it requires running multiple adjoint and direct simulations in series. The Lanczos approach can therefore still result in significant computational time, especially when the rank of the TSVD approximation needs to be increased to improve matches to observations.

\begin{figure}
\centering
\includegraphics[width=0.97\textwidth]{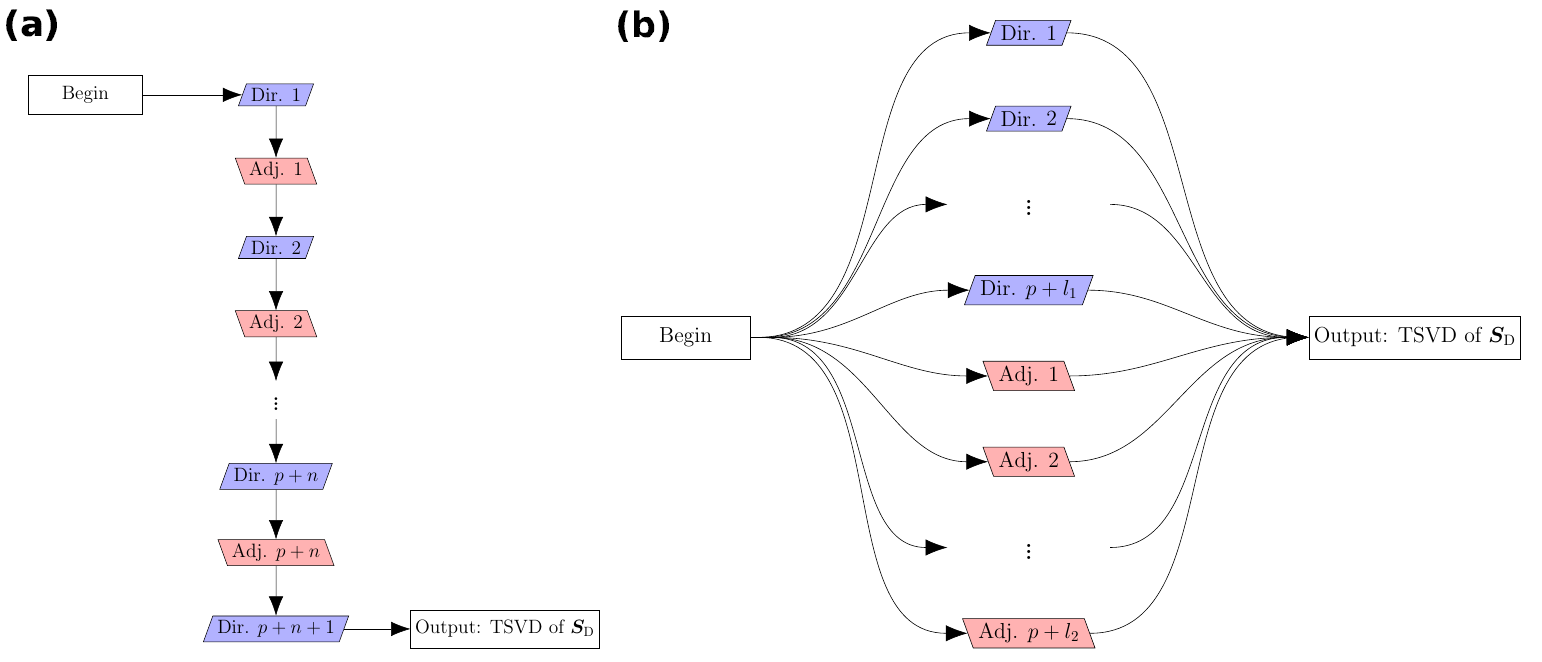}
\caption{(a) Lanczos process for a rank-$p$ TSVD of $\bm{S}_\text{D}$ requiring a series of adjoint (Adj.) and direct (Dir.) runs. (b) Randomized 1-view method running all adjoint and direct problems simultaneously in parallel to find a TSVD of $\bm{S}_\text{D}$. The integers $n$, $l_1$ and $l_2$ are oversampling factors for improved estimation.}
\label{fig:1}
\end{figure}

The present study addresses the above drawback of the TSVD-LM method proposed by \citet{tavakoli2010} by developing a new modified version of their approach which uses randomized low-rank approximation methods to evaluate the TSVD of the dimensionless sensitivity matrix. The randomized TSVD-LM approach enables simultaneous solution of the adjoint and direct simulations required at every TSVD-LM inversion iteration, which reduces computational time. The two contrasting approaches are schematically represented in Fig. \ref{fig:1}.

The Lanczos and randomized approaches are compared and contrasted in the following sections with the aim of improving the TSVD-LM approach. Section \ref{sec:background} discusses previous methods applied to speed up LM based reservoir inversions and provides a literature review motivating the use of the randomized methodology. Section \ref{sec:TLM} outlines the main aspects of the TSVD-LM method used in this study. The following section \ref{sec:TSVD} gives details for the Lanczos approach and the suggested alternative randomized TSVD methods. The methods were tested by inverting a synthetic geothermal natural state and production history model. The computational experiments are outlined in section \ref{sec:expdesign} and the inversion test results are given in section \ref{sec:results}. The results show that replacing the standard Lanczos approach with randomized methods can reduce the time spent on inversions by an order of magnitude. A summary and conclusions are given in section \ref{sec:conclusions}.

\section{Background} \label{sec:background}
\subsection{Levenberg-Marquardt Method}
In the iterative LM approach to minimizing \eqref{eq:objective}, model updates $\delta {\bm m}$ are found by solving a linear matrix equation given by
\begin{linenomath*}
\begin{equation} \label{eq:LM}
	\left[ {\bm S}^T {\bm \Gamma}^{-1}_{\text{d}} {\bm S}  + \mu {\bm R} + \gamma {\bm D}  \right]  \delta {\bm m}  =  - {\bm S}^T {\bm \Gamma}^{-1}_{\text{d}}  \left[ {\bm d} ({\bm m}) - {\bm d}_{\text{obs}}  \right]  - \mu {\bm R} \left[ {\bm m} - {\bm m}_{\text{pr}}  \right]   \, .
\end{equation} \end{linenomath*}
Here ${\bm D}$ is a positive definite matrix and $\gamma > 0$ is the adjustable LM damping parameter. Defining the observation residual vector as
\begin{linenomath*}
\begin{equation}
	{\bm r} ({\bm m}) = {\bm d} ({\bm m}) - {\bm d}_{\text{obs}}
\end{equation} \end{linenomath*}
then the sensitivity matrix ${\bm S}$ is given by
\begin{linenomath*}
\begin{equation}
\label{eq:sensitivitymatrix}
\bm{S} = \frac{\mathrm{d} [\bm{d} (\bm{m})]}{\mathrm{d} \bm{m}} = \frac{\mathrm{d} \bm{r} (\bm{m})}{\mathrm{d} \bm{m}}
= 
\begin{bmatrix}
\frac{\mathrm{d} r_1}{\mathrm{d} m_1} & \cdots & \frac{\mathrm{d} r_1}{\mathrm{d} m_{N_m}} \\
\vdots & \ddots & \vdots\\
\frac{\mathrm{d} r_{N_d}}{\mathrm{d} m_1} & \cdots & \frac{\mathrm{d} r_{N_d}}{\mathrm{d} m_{N_m}} \\
\end{bmatrix}
\, .
\end{equation} \end{linenomath*}

\subsection{Reducing Cost of Levenberg-Marquardt Updates}
The cost of solving the linear LM update equations \eqref{eq:LM} can be reduced, for example, in PEST by applying an SVD or LSQR solver. However, PEST still requires the explicit generation of the sensitivity matrix. To reduce the computational burden further, PEST provides an option called SVD-assist \citep{doherty2016} which re-parameterizes the inverse problem based on the right-singular vectors of a weighted sensitivity matrix. The SVD-assist strategy uses truncation to reduce the number of parameters, which can make inversion more manageable. This was the approach adopted by \citet{austria2015} and \citet{colina2013} for improving geothermal models, with sensitivity values evaluated using finite differencing. A drawback of the SVD-assist approach, as currently implemented in PEST, is that it requires forming the sensitivity matrix for the full parameterization at least once.

Since evaluating the sensitivity matrix is often too expensive, the quasi-Newton limited memory BFGS method which only requires the gradient of the objective function for every model update has been preferred for modeling petroleum reservoirs \citep{gao2006improved,zhang2002optimization}. The gradient can be found efficiently using adjoint codes at a cost not exceeding that of one model simulation, independent of the number of parameters.

\citet{tavakoli2010} demonstrated that the LM approach can be made computationally more efficient than the limited memory BFGS method by basing model updates on the principal right-singular vectors of a dimensionless sensitivity matrix. Their method avoids explicitly forming the sensitivity matrix by applying the Golub-Kahan-Lanczos bidiagonalization method \citep{golub2013,vogel1994} to find a TSVD of the dimensionless sensitivity matrix. Each iteration of the Lanczos method only requires the evaluation of the products of ${\bm S}$ times a vector and ${\bm S}^T$ times a vector. The former can be evaluated efficiently using the direct method \citep{oliver2008,rodrigues2006adjoint} while the latter can be found at a similar cost using the adjoint method \citep{oliver2008,rodrigues2006adjoint}. Moreover, the TSVD method gave better parameters with a lower model mismatch (Eq. \eqref{eq:modelmismatch}) than the BFGS method \citep{tavakoli2010}.

Expanding on their work further, \citet{tavakoli2011} applied their TSVD method to finding ensembles of models conditioned on data using the randomized maximum likelihood method (RML) \citep{kitanidis1995,oliver1996RMLshortnote,oliver1996conditioning}. The RML method was developed as a relatively inexpensive (compared to Markov chain Monte Carlo) approximate posterior sampler for reservoir modeling. Though RML has only been proven to sample correctly for linear problems \citep{oliver1996RMLshortnote,oliver2008}, results from numerical experiments indicate that the method is applicable to nonlinear forward problems \citep{emerick2013investigation,gao2006punqs3,liu_2003,shirangi2014,shirangi2016,tavakoli2011,zafari2007}. 

RML requires solving a costly inverse problem for every ensemble member. However, \citet{tavakoli2011} showed that using the TSVD approach for RML sampling can reduce the computational cost of approximating the posterior. Their method called SVD-EnRML achieved this by evaluating at every LM iteration a TSVD of the dimensionless sensitivity matrix for only a base set of model parameters and using this base TSVD to propose candidate parameter updates for all ensemble members.

\citet{tavakoli2010,tavakoli2011} applied the TSVD-LM method to inferring permeability distributions of two-dimensional synthetic petroleum reservoirs using production data. \citet{shirangi2014} later on used the method for history matching permeability and porosity values in three-dimensional synthetic petroleum reservoir models. \citet{shirangi2014} enhanced the SVD-EnRML algorithm by introducing an ensemble-based regularization method which can be applied in cases where an ensemble of models can be drawn from the prior though a prior parameter covariance is too large to be generated and stored. \citet{shirangi2016} made further improvements to the SVD-EnRML method as well as looking into the difference between using LM parameter updates and those obtained with the Gauss-Newton (GN) method.

Iterative Krylov linear solvers such as CG or LSQR can also be applied to solve the linear LM update equations. Like the TSVD approach of \citet{tavakoli2010} these methods only require the sensitivity matrix times a vector and the sensitivity matrix transposed times a vector at every Krylov iteration. Since approximate solutions will often suffice, the Krylov linear solvers can be terminated early when a maximum number of iterations has been reached or the approximate solution satisfies a predefined tolerance. The CG and LSQR methods may therefore enable inversion of a single reservoir model with an efficiency similar to the TSVD-LM method of \citet{tavakoli2010}. 

However, the LM updates depend on the LM damping parameter $\gamma$. As a result when a candidate update fails and the damping parameter is updated then straightforward implementations using CG or LSQR require a repeated solution of the LM update equations accompanied by additional adjoint and direct solves. The TSVD-LM method proposed by \citet{tavakoli2010} does not have this drawback since the TSVD of the dimensionless sensitivity matrix can be re-used, independent of the LM damping factor. 

Considering all the above points, straightforward CG or LSQR implementations might be expected to result in more adjoint/direct solves and slower implementations than the TSVD-LM method. However, this is likely to be problem specific as CG or LSQR solvers will be especially effective when the LM update equations are well-conditioned or can be made so using a good choice of preconditioner.

\citet{lin2016} recently developed a LM method which uses the LSQR method in a more sophisticated way. Their method finds a Krylov subspace using Golub-Kahan-Lanczos bidiagonalization, independent of the LM damping parameter. The subspace is then used to solve for model updates using the LSQR method. Though \citet{lin2016} implemented their method for problems where the full sensitivity matrix can be found without too much expense, their Levenberg variant, with the LM damping matrix ${\bm D}$ chosen as the identity matrix, only requires ${\bm S}$ and ${\bm S}^T$ times vectors during the bidiagonalization procedure. By evaluating ${\bm S}$ and ${\bm S}^T$ times vectors using adjoint and direct methods, the LSQR-LM approach of \citet{lin2016} may therefore have similar characteristics and efficiency as the TSVD-LM method of \citet{tavakoli2010}.

However, a drawback of the Golub-Kahan-Lanczos bidiagonalization method applied by \citet{lin2016,shirangi2014,shirangi2016,tavakoli2010,tavakoli2011} is that it is inherently serial. For a rank-$p$ SVD approximation the method requires $p$ or more iterations, each requiring an adjoint run followed by a direct solve. The method could be made much more efficient if the sensitivity matrix to vector products could all be evaluated simultaneously in parallel. 

Randomzied matrix approximation methods are considered here because they allow simultaneous evaluation of all sensitivity matrix to vector products.

\subsection{Randomized Alternatives}
Randomized methods \citep{erichson2016,halko2011structure,martinson2016rsvd,rokhlin2009,woolfe2008} have been studied extensively in recent years as a means of low-rank SVD matrix approximation. Various randomized SVD algorithms have been developed that promise to alleviate the lack of parallelizability of classical SVD methods, while returning a similar level of accuracy as the standard Lanczos method \citep{halko2011structure}. 

Instead of dealing with matrix-vector products randomized SVD methods use a small number of matrix-matrix products which means that the methods lend themselves more readily to parallel and high performance computing \citep{erichson2016,gu2015,halko2011structure,martinson2016rsvd,voronin2015}. Randomized methods allow for faster SVD approximations by reducing the number of times information contained in the matrix of interest is accessed. It may suffice to access or view the matrix of interest twice \citep{halko2011structure} or just once \citep{halko2011structure,martinson2016rsvd,tropp2016,woolfe2008}, which is especially advantageous for matrices too big to fit into core memory. Randomized SVD methods have therefore been applied to principal component analysis (PCA) of large data sets \citep{halko2011pca}.

In the modeling context, randomized SVD methods have been applied to PCA of covariance matrices \citep{dehdari2012,lee2014,saibaba2015}, inversion \citep{lee2014,xiang2013,lee2016} and uncertainty quantification \citep{buithanh2012,cui2016,cui2014,isaac2015}. \citet{buithanh2012} and \citet{isaac2015} used a matrix-free Newton-CG method to invert a global seismic model and an Antarctic ice sheet model, respectively. Both studies used randomized methods for low-rank approximation of posterior Hessian matrices for over a million parameters. \citet{buithanh2012} and \citet{isaac2015} used adjoint and direct methods to evaluate the action of the Hessian on matrices, required for the randomized approximation.

\citet{lee2014,lee2016} and \citet{saibaba2015} used the quasi-linear geostatiscal approach (GA) \citep{kitanidis1995} to invert hydrological models. \citet{lee2014} applied randomized PCA on the prior covariance matrix to reduce the computational cost of inversion. \citet{saibaba2015}, however, used a randomized approach to approximate the posterior covariance matrix. \citet{lee2016} pointed out that GA can be improved further by using randomized low-rank matrix approximations to precondition the linear GA update equations.
\citet{lin2017} demonstrated, similarly, that the scalability of GA can be substantially improved by using randomized sketching to reduce the effective size of the observation space. 

\citet{xiang2013} presented a method for solving Tikhonov regularized inverse problems by finding a TSVD of the sensitivity matrix using randomized methods and solving the GN update equations using a pseudoinverse. They demonstrated the method on linear forward problems though the approach may be extended to nonlinear forward problems.

Stochastic gradient methods have been developed for history-matching of reservoir models when adjoint code is not available \citep{fonseca2016,gao2007,li2011}. These randomized gradient methods inherit the drawback of steepest descent methods, i.e., a relatively slow convergence rate compared with LM and other Newton-like methods, with the added drawback of using approximate gradients.

\subsection{Summary of Approach Developed Here}
The present study looks at applying randomized SVD methods presented in \citet{halko2011structure} and \citet{tropp2016} along with adjoint and direct code to speed up the TSVD-LM scheme discussed in \citet{tavakoli2010}. 

Using a randomized method to form a TSVD of the dimensionless sensitivity matrix requires the same type and similar number of sensitivity matrix to vector products (adjoint and direct problems) as the Lanczos method. However, computational savings are possible since the sensitivity products (adjoint and direct problems) can be made independent of each other and evaluated simultaneously.

As previously used by \citet{isaac2015} and \citet{buithanh2012}, we use adjoint and direct methods for efficient evaluation of the action of the sensitivity matrix (or its transpose) on matrices. Previous work by \citet{lee2014,lin2017,saibaba2015,buithanh2012,isaac2015} used randomized sketching either when initializing inversions or after running inversions. This is unlike the present study which applies new randomized sketching at every inversion (LM) iteration. The randomized preconditioner proposed by \citet{lee2016}, to improve the scalability of GA, may likewise be formed at every inversion iteration. However, \citet{lee2016} used finite-differencing to evaluate the sensitivity matrix products needed by their method.

Unlike previous studies, the present one applies randomized TSVD methods, adjoint code and direct code at every TSVD-LM iteration to speed up inversion. Furthermore, this is the first study on inversion of reservoir models to apply a randomized method recently developed by \citet{tropp2016}. As detailed in section \ref{sec:1view}, the method of \citet{tropp2016} enables solving all the adjoint problems in parallel and concurrently with solving the direct problems. Section \ref{sec:subreuse} presents new randomized TSVD variants based on subspace ideas proposed by \citet{vogel1993}.

\section{Truncated Levenberg-Marquardt} \label{sec:TLM}
\subsection{Model Updates}
As mentioned above, forming \eqref{eq:LM} and solving it directly may be computationally prohibitive. Following the work of \citet{tavakoli2010}, transformed model parameters are found according to
\begin{linenomath*}
\begin{equation} \label{eq:partrans}
	\tilde{\bm m} = {\bm L}^{-1} \left[ {\bm m} - {\bm m}_\text{pr}  \right]  \, .
\end{equation} \end{linenomath*}
Here ${\bm L} = ({\bm R}^{-1})^{1/2}$ is the square root matrix of ${\bm R}^{-1}$. \citet{tavakoli2010} chose ${\bm R}^{-1}$ as a prior covariance ${\bm \Gamma}_\text{m}$ and suggested using a Cholesky decomposition such that ${\bm R}^{-1} = {\bm L} {\bm L}^T$ where ${\bm L}$ is a lower triangular matrix. However, for a large number of parameters creating and storing a Cholesky matrix ${\bm L}$ will be computationally expensive. Alternatively, if the eigenvalues of ${\bm R}^{-1}$ decay rapidly enough, the matrix ${\bm L}$ can be approximated using a truncated eigendecomposition of ${\bm R}^{-1}$. The truncated eigendecomposition or principal components can be found efficiently using randomized methods \citep{dehdari2012,halko2011pca,halko2011structure,lee2014,martinson2016rsvd,saibaba2015} when ${\bm R}^{-1}$ times a thin matrix can be evaluated efficiently. When forming a prior covariance is too expensive then regularization can also be based on a low-rank matrix $({\bm R}^{-1})^{1/2}$ generated by an ensemble of models sampled from the prior \citep{shirangi2014}. The Cholesky approach was used in the present study because the size of the inverse problem tested in section \ref{sec:results} is small enough so that the Cholesky approach is feasible. 

Here we will use the common assumption that the observation covariance is diagonal so it can be written as ${\bm \Gamma}_\text{d} = {\bm \Gamma}_\text{d}^{1/2} {\bm \Gamma}_\text{d}^{T/2}$ where ${\bm \Gamma}_\text{d}^{1/2}$ is a diagonal matrix with each diagonal element equal to the square root of the corresponding element of ${\bm \Gamma}_\text{d}$. Then using the transformed parameters \eqref{eq:partrans} the objective function can be written as \citep{shirangi2011}
\begin{linenomath*}
\begin{equation}
	\Phi (\tilde{\bm m}) = \left[ {\bm \Gamma}_\text{d}^{-1/2}  {\bm r} (\tilde{\bm m})   \right]^T  \left[ {\bm \Gamma}_\text{d}^{-1/2}  {\bm r} (\tilde{\bm m})   \right]  +  \mu \tilde{\bm m}^T \tilde{\bm m}  \, .
\end{equation} \end{linenomath*}

For the transformed parameters the LM update equations are
\begin{linenomath*}
\begin{equation} \label{eq:LMtrans1}
	\left[ {\bm L}^T {\bm S}^T {\bm \Gamma}^{-1}_{\text{d}} {\bm S} {\bm L}  + \mu {\bm I} + \gamma \tilde{\bm D}  \right]  \delta \tilde{\bm m}  =  - {\bm L}^T {\bm S}^T {\bm \Gamma}^{-1}_{\text{d}}   {\bm r} (\tilde{\bm m})  - \mu \tilde{\bm m}   \, ,
\end{equation} \end{linenomath*}
Choosing the LM damping matrix $\tilde{\bm D}$ as the identity and introducing the dimensionless sensitivity matrix ${\bm S}_\text{D} = {\bm \Gamma}_\text{d}^{-1/2} {\bm S} {\bm L}$ \citep{tavakoli2010,zhang2002} then Eq. \eqref{eq:LMtrans1} can be written as
\begin{linenomath*}
\begin{equation} \label{eq:LMtrans3}
	\left[ {\bm S}_\text{D}^T {\bm S}_\text{D}  + (\mu + \gamma ) {\bm I}\right]  \delta \tilde{\bm m}  =  - {\bm S}_\text{D}^T {\bm \Gamma}^{-1/2}_{\text{d}}   {\bm r} (\tilde{\bm m})  - \mu \tilde{\bm m}   \, .
\end{equation} \end{linenomath*}

The SVD factorization of ${\bm S}_\text{D}$ is
\begin{linenomath*}
\begin{equation}
	{\bm S}_\text{D} = {\bm U} {\bm \Lambda} {\bm V}^T = \sum_{i=1}^N \lambda_i {\bm u}_i {\bm v}_i^T   \, ,
\end{equation} \end{linenomath*}
where $N = \min(N_d , N_m)$, the matrix ${\bm U} = [ {\bm u}_1 \,  {\bm u}_2 \, \dots \, {\bm u}_N ]$ is orthogonal, the matrix ${\bm V} = [ {\bm v}_1 \,  {\bm v}_2 \, \dots \, {\bm v}_N ]$ is orthogonal and ${\bm \Lambda} = \text{diag} [ \lambda_1 \, ,  \lambda_2 \, , \, \dots \, , \lambda_N ]$, with $\lambda_1 \geq \lambda_2  \geq \dots \geq \lambda_N$. $\lambda_i$ is the $i$th singular value while ${\bm u}_i$ and ${\bm v}_i$ are its corresponding left- and right-singular vectors, respectively. 

Truncating the SVD and retaining the singular triplets $(\lambda_i, {\bm u}_i, {\bm v}_i)$ for the $p$ largest singular values we can make the following rank-$p$ approximation of ${\bm S}_\text{D}$
\begin{linenomath*}
\begin{equation} \label{eq:TSVD}
	{\bm S}_\text{D} \approx \sum_{i=1}^p \lambda_i {\bm u}_i {\bm v}_i^T  = {\bm U}_p {\bm \Lambda}_p {\bm V}_p^T  \, .
\end{equation} \end{linenomath*}
Using the TSVD \eqref{eq:TSVD}, \citet{tavakoli2010} showed that an approximate solution to \eqref{eq:LMtrans3} can be found by basing model updates on the principle right-singular vectors. \citet{tavakoli2010} looked at the case with the regularization weight $\mu = 1$. For the more general case with a variable regularization weight their update equations can be written as
\begin{linenomath*}
\begin{equation} \label{eq:TSVDLM}
	\delta \tilde{\bm m} = \sum_{i=1}^p \alpha_i {\bm v}_i  \, ,
\end{equation} \end{linenomath*}
where
\begin{linenomath*}
\begin{equation} \label{eq:alpha}
	\alpha_i = - \frac{1}{\mu + \gamma + \lambda_i^2}  \left[  \mu {\bm v}_i^T \tilde{\bm m}  + \lambda_i {\bm u}_i^T {\bm \Gamma}^{-1/2}_{\text{d}}   {\bm r} ({\bm m}) \right]  \, .
\end{equation} \end{linenomath*}
The model update in the original parameter space is $\delta {\bm m} = {\bm L} \delta \tilde{\bm m}$.

Various strategies are available for determining $\mu$ \citep{aster2005,haber1997thesis,vogel2002}. The regularization weight $\mu$ can, for instance, be gradually decreased for a sequence of separate inverse problems \citep{haber1997thesis,haber2000inv}. Alternatively, to reduce computational cost, a single inversion can be run with the weight varied at every iteration using a local linear approximation of the objective function \citep{doherty2003,doherty2016,haber1997thesis,haber2000gcv}. Though adaptive adjustment of the regularization weight is a helpful option we do not consider it further here and use a fixed $\mu$. The LM damping factor $\gamma$ and the number of retained singular triplets are instead used as the primary means of regulating model updates at each LM iteration, (see following section).

To regularize the inverse problem further we also prescribe upper and lower bounds for the model parameters. This is especially helpful for geothermal problems since extreme parameter values commonly result in slow model runs and may result in run failure.

If the singular triplets are exact then \eqref{eq:TSVD} is the optimal rank-$p$ approximation of ${\bm S}_\text{D}$ in terms of the spectral and Frobenius norms \citep{golub2013}. However, for large problems the computational cost of finding precise singular triplets is typically not worth the effort. Section \ref{sec:TSVD} discusses methods for low-cost approximation of the singular triplets of ${\bm S}_\text{D}$.

\subsection{The Levenberg-Marquardt Damping Factor} \label{sec:LMdamping}
The LM damping factor $\gamma$ controls whether the model updates are closer to steepest descent ($\gamma \rightarrow \infty$) or a GN update ($\gamma \rightarrow 0$). It is usual to start off with a large value for $\gamma$ at early iterations for more globally robust steepest-descent-like updates and to gradually reduce $\gamma$ as the LM method proceeds to enjoy the faster asymptotic convergence rate of the GN method.

The choice of $\gamma$ not only controls the inversion procedure, but also regularizes each model update. Applying GN updates, or using small values of $\gamma$, can introduce unwanted inversion artifacts since it admits large contributions from the singular vectors associated with small singular values, which introduce high frequency components to the model parameter updates. 

\citet{gavalas1976} observed for history matching of a one-dimensional flow model that using a larger LM damping factor resulted in spatially smoother porosity and permeability values. \citet{abacioglu2001} found that history matching two- and three-dimensional reservoir models using the GN method ($\gamma = 0$) or a LM approach where $\gamma$ is small early on could result in model roughness which was difficult to remove. Furthermore, they found that LM updates with small $\gamma$ can result in convergence to a large final objective function value \citep{abacioglu2001}. To combat these issues they recommended initializing with a large $\gamma$; they also tried using an adjustable regularization weight $\mu$ with the aim of keeping the contribution of the model mismatch term similar to the observation mismatch. 

Similarly, \citet{gao2006improved} found when applying a limited memory BFGS method to history matching that adaptively re-weighting the observation mismatch term can prevent bad observation matches because of model over adjustment at early iterations. 

\citet{li2003} compared the GN method against the LM method for history matching a two-dimensional reservoir. Their results showed that the GN method resulted in a poorer match and a rougher permeability distribution than the LM approach, despite the GN method using twice as many iterations.

Reducing the number of right-singular vectors used for the model updates in \eqref{eq:TSVDLM}, by truncating away the singular vectors associated with small singular values, can also help to regularize the inverse problem. Applying a TSVD-GN method may therefore avoid some of the pitfalls of full GN updates. With this in mind, \citet{shirangi2016} compared LM and GN variants of their SVD-EnRML approach for finding reservoir models conditioned on production history data. Though the GN variant of SVD-EnRML can be improved by appropriate truncation, \citet{shirangi2016} concluded that the inclusion of $\gamma$ makes the LM approach more robust to the number of singular triplets that are used. This makes TSVD-LM preferable to TSVD-GN since the strategy for determining the number of retained singular triplets requires input from the modeler, which is likely to result in suboptimal truncation.

\begin{algorithm}[t!]
  \caption{TSVD Levenberg-Marquardt}\label{alg:TSVDLM}
  \textbf{INPUT:} {Initial parameter guess ${\bm m}_\text{pr} \in {\mathbb{R}}^{N_m}$, $iter_{\text{max}}$ and $\varepsilon_m$.}\\
  \textbf{RETURNS:} {Model parameters ${\bm m}$.}
  \begin{algorithmic}[1]	
    \State ${\bm m} = {\bm m}_{\text{pr}}$; {\ttfamily converged} = \textbf{False}; $iter = 1$; $\gamma = 10^6$.
    \State Evaluate ${\bm r} ({\bm m})$ by running a nonlinear simulation.
    \State Estimate a TSVD of ${\bm S}_{\text{D}} ({\bm m})$.
    \State $\tilde{\bm m} = {\bm L}^{-1} \left[ {\bm m} - {\bm m}_{\text{pr}} \right]$.
    \While {(\textbf{not}  {\ttfamily converged}) \textbf{and} ($iter \leq iter_{\text{max}}$)}
    	\State Find model update $\delta \tilde{\bm m}$ using Eq. \eqref{eq:TSVDLM}.
        \State ${\bm m}_{\text{temp}} = {\bm m} + {\bm L} \delta \tilde{\bm m}$.
        \State Truncate ${\bm m}_{\text{temp}}$ to be within the parameter bounds.
        \If {( $\left\lVert{\bm m}_\text{temp} - {\bm m}\right\rVert \leq \varepsilon_m [ \left\lVert{\bm m}\right\rVert + \varepsilon_m ]$ )}
        	\State {\ttfamily converged} = \textbf{True}.
        \Else 
        	\State Evaluate ${\bm r} ({\bm m}_{\text{temp}})$ by running a nonlinear simulation.
            \State $\tilde{\bm m}_{\text{temp}} = {\bm L}^{-1} \left[ {\bm m}_{\text{temp}} - {\bm m}_{\text{pr}} \right]$.
            \If {( $\Phi ({\bm m}_{\text{temp}}) < \Phi ({\bm m}) $)}
            	\State $iter = iter + 1$; ${\bm m} = {\bm m}_{\text{temp}}$; $\tilde{\bm m} = \tilde{\bm m}_{\text{temp}}$.
                \State $\gamma = \gamma / 10$.
                \State Estimate a TSVD of ${\bm S}_{\text{D}} ({\bm m})$.
            \Else
            	\State $\gamma = \max (10 \gamma , 100)$.
            \EndIf
        \EndIf
    \EndWhile
  \end{algorithmic}
\end{algorithm}

The above discussion suggests that the LM approach is likely to result in better matches to observations than the GN approach. Furthermore, the LM approach should generally result in models which honor prior knowledge more closely. Therefore, we conclude that including the damping factor $\gamma$ makes the LM approach more robust than and preferable to the GN approach for inverting reservoir models. The TSVD-LM method is therefore adopted here instead of using a TSVD-GN method.

The LM damping parameter $\gamma$ can be chosen to have a large value at early iterations and gradually reduced at subsequent iterations to improve the model resolution. Algorithm \ref{alg:TSVDLM} outlines the TSVD-LM method used here. It is based on the same sequential way of varying $\gamma$ as \citet{shirangi2014} and \citet{shirangi2016}. If a candidate model update ${\bm m}_\text{temp} = {\bm m} + \delta {\bm m}$ results in a lower value for the objective function, then we set ${\bm m} = {\bm m}_\text{temp}$ and move to the next iteration with $\gamma = \gamma / 10$. Otherwise the damping factor is increased according to $\gamma = \max (10 \gamma , 100)$. The factor is increased to at least $100$ when a model update fails in an attempt to reduce the number of forward simulations, since a small value of $\gamma$ will not have much of an impact on the denominator in \eqref{eq:alpha}. A good initial $\gamma$ is problem specific and selecting a damping factor that is too large can lead to slow convergence. Here we use an initial $\gamma = 10^6$ as suggested by \citet{shirangi2014} and \citet{shirangi2016}.

The sequential method adopted here for varying $\gamma$ was chosen for its simplicity. However, because it is sequential this approach can result in slow progress if the proposed LM damping factors lead to multiple unsuccessful model updates or small reductions in the objective function. Another approach, not used here, is to find simultaneously in parallel model updates and corresponding model outputs for multiple damping factors \citep{doherty2016,lin2016} which may reduce computational time and improve convergence performance. The TSVD-LM method, like the method proposed by \citet{lin2016}, is well suited to this parallel approach since no additional adjoint or direct solves are needed when $\gamma$ is varied during a LM iteration. However, straightforward CG or LSQR implementations result in additional adjoint and direct solves for every $\gamma$.

\section{TSVD of the Dimensionless Sensitivity Matrix} \label{sec:TSVD}
\subsection{Truncated SVD Using Lanczos Method}
For finding the TSVD of ${\bm S}_\text{D}$ \citet{shirangi2014,shirangi2016,tavakoli2010,tavakoli2011} applied Lanczos iteration \citep{golub2013,vogel1994}. Algorithm \ref{alg:lanczos} outlines a basic Lanczos algorithm based on the one given by \citet{vogel1994} for finding a TSVD of a matrix. \citet{vogel1994} applied Lanczos iteration to find an approximate solution to an example linear inverse problem.

\begin{algorithm}
  \caption{Lanczos Method for TSVD}\label{alg:lanczos}
  \textbf{INPUT:} {Matrix ${\bm A} \in {\mathbb{R}}^{n_r \times n_c}$, integer $p > 0$ and convergence tolerance $\varepsilon_\text{sv}$.}\\
  \textbf{RETURNS:} {Approximate rank-$p$ SVD, ${\bm U}_p {\bm \Lambda}_p {\bm V}_p^T$, of ${\bm A}$.}
  \begin{algorithmic}[1]	
    \State Generate unit vector ${\bm q}_1 \in {\mathbb{R}}^{n_c}$.
    \State Compute ${\bm y} = {\bm A} {\bm q}^1$ ; $\alpha_1 = \left\lVert{\bm y}\right\rVert$ ; ${\bm p}^1 = {\bm y}/\alpha_1$.   \label{algline:lancdir1}
    \State Define ${\bm Q}^j = \left[ {\bm q}^1, \dots , {\bm q}^j \right]$ and ${\bm P}^j = \left[ {\bm p}^1, \dots , {\bm p}^j \right]$.
    \State Set $j=1$ , {\ttfamily converged} = \textbf{False} and $\lambda^j_i = 0$, for $i=1, \, 2, \, \dots , \, p$.
    \While {(\textbf{not}  {\ttfamily converged})}
        \State ${\bm w} = {\bm A}^T {\bm p}^j - \alpha_j {\bm q}^j$.    \label{algline:lancadjj}
        \State Reorthogonalize: ${\bm w} = {\bm w} - {\bm Q}^j \left( \left[ {\bm Q}^j  \right]^T {\bm w} \right)$.
        \State $\beta_j = \left\lVert{\bm w}\right\rVert$; ${\bm q}^{j+1} = {\bm w}/\beta_j$.
        \State ${\bm y} = {\bm A} {\bm q}^{j+1} - \beta_j {\bm p}^j$.    \label{algline:lancdirj}
        \State Reorthogonalize: ${\bm y} = {\bm y} - {\bm P}^j \left( \left[ {\bm P}^j  \right]^T {\bm y} \right)$.
        \State $\alpha_{j+1} = \left\lVert{\bm y}\right\rVert$; ${\bm p}^{j+1} = {\bm y}/\alpha_{j+1}$.
        \State Define the bidiagonal matrix $\tilde{\bm T}^{j+1} = \begin{bmatrix}
        \alpha_1 & \beta_1 & & \\
         & \alpha_2 & \ddots & \\
         & & \ddots & \beta_{j} \\
         & & & \alpha_{j+1}
        \end{bmatrix}$.
        \If {($j \geq p$)}
        	\State Evaluate the SVD of $\tilde{\bm T}^{j+1} = \tilde{\bm U} {\bm \Lambda}_{j+1} \tilde{\bm V}^T$ and truncate.
            \If {(Eq. \eqref{eq:lancConv} holds)}
            	\State {\ttfamily converged} = \textbf{True}.
            \EndIf
        \EndIf
        \State $j=j+1$.
    \EndWhile
    \State ${\bm U}_p = {\bm P}^j \tilde{\bm U}_p$ and ${\bm V}_p = {\bm Q}^j \tilde{\bm V}_p$.
  \end{algorithmic}
\end{algorithm}

The main cost involved in finding the TSVD of ${\bm S}_\text{D}$ using the Lanczos approach is from the matrix vector multiplications in lines \ref{algline:lancdir1}, \ref{algline:lancadjj} and \ref{algline:lancdirj} of Algorithm \ref{alg:lanczos}. For evaluating ${\bm S}_\text{D} {\bm q}^j$ we begin by finding $\tilde{\bm q} = {\bm L} {\bm q}^j$ followed by evaluating $\hat{\bm q} = {\bm S} \tilde{\bm q}$, which can be found efficiently using the direct method. Finally, ${\bm S}_\text{D} {\bm q}^j = {\bm \Gamma}_\text{d}^{-1/2} \hat{\bm q}$. 
When calculating ${\bm S}_\text{D}^T {\bm p}^j$ we first find $\tilde{\bm p} = {\bm \Gamma}_\text{d}^{-T/2} {\bm p^j}$ followed by $\hat{\bm p} = {\bm S}^T \tilde{\bm p}$, which is found efficiently using the adjoint method, and then ${\bm S}_\text{D}^T {\bm p}^j = {\bm L}^T \hat{\bm p}$. The adjoint and direct methods for finding ${\bm S}^T$ and ${\bm S}$ times vectors are outlined in Appendices A and B.

To approximate the first $p$ singular triplets of ${\bm S}_\text{D}$ the Lanczos procedure can be halted when the number of Lanczos iterations $j$ exceed or equal $p$ and
\begin{linenomath*}
\begin{equation} \label{eq:lancConv}
	\frac{\lvert \lambda_i^{j+1} - \lambda_i^{j} \rvert}{\lambda_i^{j+1}} \leq \varepsilon_\text{sv} \quad , \quad \text{for} \quad i = 1 \, , 2 \, , \, \dots \, , p
\end{equation} \end{linenomath*}
\citep{vogel1994}. To achieve the desired precision \eqref{eq:lancConv} requires $p + n$ Lanczos iterations. Tavakoli and Reynolds [2010] and Shirangi [2011] reported that $n$ varied between 3 and 8 when using $\varepsilon_\text{sv} = 10^{-6}$ or $10^{-5}$ for their problems. The Lanczos procedure requires a series of $p+n+1$ direct and $p+n$ adjoint solves, which becomes costly when retaining a large number of singular triplets.

For reducing the computational cost of the TSVD method \citet{tavakoli2010} suggested using a small number of singular triplets $p$ at early iterations and gradually increasing $p$ as the inversion proceeds. This can for instance be achieved by increasing $p$ between iterations by some fixed integer. Algorithm \ref{alg:lanczos} is presented for this case where $p$ can be given as input.

Another approach introduced by \citet{tavakoli2010} and subsequently used by \citet{shirangi2014,shirangi2016,tavakoli2011} is to truncate based on the ratio of the largest singular value to the retained singular values. Then the truncation $p$ can be chosen as the smallest value such that
\begin{linenomath*}
\begin{equation} \label{eq:svcut}
	\frac{\lambda_p}{\lambda_1} \leq \text{sv-cut}
\end{equation} \end{linenomath*}
\citep{shirangi2014,shirangi2016,tavakoli2011} and sv-cut can be gradually decreased between LM iterations to include more singular triplets. 

As $p$ is increased and $\gamma$ is lowered during the inversion higher frequency components are gradually allowed greater influence on model updates to introduce finer spatial details in the model which may be required for a good match to observations. The results of \citet{shirangi2014,shirangi2016,tavakoli2010,tavakoli2011} indicate that regulating the retained number of singular triplets by lowering sv-cut is a sound procedure. 

Truncation control using sv-cut is convenient when using the Lanczos method since approximate singular values can be evaluated at every Lanczos iteration. A check can then be made at each Lanczos iteration to see whether \eqref{eq:svcut} and \eqref{eq:lancConv} are fulfilled.

However, the sv-cut approach is not used for the randomized TSVD methods because the randomized methods considered in the present study are non-iterative. Instead, when using randomized methods we gradually increase the number of retained singular triplets $p$ by a fixed value between LM iterations. The suitability of this strategy is investigated in section \ref{sec:resultstrunccontrol} by comparing it with the sv-cut approach.

\subsection{Randomized TSVD Methods}
A drawback of the Lanczos approach is that it necessitates at least $p$ computationally expensive iterations. Each iteration requires one direct solve and one adjoint solve. Evaluating the TSVD of ${\bm S}_\text{D}$ by this serial Lanczos approach is therefore very time-consuming. 

Randomized methods are a promising alternative for estimating the rank-$p$ approximation of ${\bm S}_\text{D}$ since they are parallelizable. Variants of the randomized approach are discussed below.

\subsubsection{Randomized 2-View Method}
Algorithm \ref{alg:2view} provides an elementary randomized method for evaluating an approximate TSVD of a matrix \citep{halko2011structure,martinson2016rsvd}. Some sources call this algorithm the basic randomized SVD algorithm \citep{gu2015,martinson2016rsvd}. Here the algorithm is called the randomized 2-view method since it only requires viewing or accessing the matrix of interest twice, once when forming the sample matrix ${\bm Y}$ in line \ref{algline:2viewYmat} and a second time for creating the small matrix ${\bm B}^T$ in line \ref{algline:2viewBmat}.

\begin{algorithm}
  \caption{Randomized 2-view method for TSVD}\label{alg:2view}
  \textbf{INPUT:} {Matrix ${\bm A} \in {\mathbb{R}}^{n_r \times n_c}$ ($n_r \geq n_c$), integers $p > 0$ and $l \geq 0$.}\\
  \textbf{RETURNS:} {Approximate rank-$p$ SVD, ${\bm U}_p {\bm \Lambda}_p {\bm V}_p^T$, of ${\bm A}$.}
  \begin{algorithmic}[1]	
    \State Generate a Gaussian random matrix ${\bm \Omega} \in {\mathbb{R}}^{n_c \times (p+l)}$.
    \State Form the matrix ${\bm Y} = {\bm A} {\bm \Omega}.$ \Comment{${\bm Y} \in {\mathbb{R}}^{n_r \times (p+l)}$}   \label{algline:2viewYmat}
    \State Find an orthonormal matrix ${\bm Q}\in {\mathbb{R}}^{n_r \times (p+l)}$, using QR factorization, such that ${\bm Y}={\bm Q} \tilde{\bm R}$.
    \State Evaluate the matrix ${\bm B}^T={\bm A}^T {\bm Q}$. \Comment{${\bm B}^T \in {\mathbb{R}}^{(p+l) \times n_c}$}    \label{algline:2viewBmat}
    \State Calculate the SVD of the relatively small matrix ${\bm B}^T ={\bm V}_{p+l} {\bm \Lambda}_{p+l} \hat{\bm U}_{p+l}^T$ and truncate.
    \State Form the matrix ${\bm U}_p = {\bm Q} \hat{\bm U}_p$.
  \end{algorithmic}
\end{algorithm}

The first three steps in Algorithm \ref{alg:2view} find an approximate range of the matrix ${\bm A}$ by finding the action of ${\bm A}$ on a random sampling matrix. Here the elements of the random sampling matrix are drawn from a standard Gaussian distribution, though other choices can be made \citep{halko2011structure,szlam2014,tropp2016}. The oversampling parameter $l$ determines the number of extra columns for the random matrix. The robustness of the method improves as $l$ is increased. Here we use $l = 10$ as it often works well \citep{martinson2016rsvd}.

Algorithm \ref{alg:2view} can be applied to ${\bm S}_\text{D}$ or its transpose, depending on the dimensions of ${\bm S}_\text{D}$, to estimate the TSVD. For reservoir simulations the main computational expense is from steps \ref{algline:2viewYmat} and \ref{algline:2viewBmat} which require evaluating ${\bm S}_\text{D}$ and ${\bm S}_\text{D}^T$ times thin matrices. The 2-view method requires the action of ${\bm S}_\text{D}$ and ${\bm S}_\text{D}^T$ on the columns of the matrices ${\bm \Omega}$ or ${\bm Q}$. These matrix products can be found using $p+l$ direct solves and $p+l$ adjoint solves.

The total number of adjoint and direct solves is nearly the same as the Lanczos approach when $l = n$. However, unlike for the Lanczos method, the $p+l$ direct solves are essentially independent of each other, and likewise for the $p+l$ adjoint solves. This produces the possibility of greatly speeding up the low-rank approximation of ${\bm S}_\text{D}$ by solving all the $p+l$ direct problems simultaneously in parallel and similarly for the adjoint problems.

Evaluating ${\bm S}_\text{D}$ times a matrix ${\bm H}$ with $p+l$ columns proceeds similarly to evaluating ${\bm S}_\text{D}$ times a vector. Here we first evaluate $\tilde{\bm H} = {\bm L} {\bm H}$ and then evaluate ${\bm S} \tilde{\bm H}$, which can be found using the direct method (see Appendix A), and then ${\bm S}_\text{D} {\bm H} = {\bm \Gamma}_\text{d}^{-1/2} \left[ {\bm S} \tilde{\bm H} \right]$. 

Here most of the computational savings are related to the evaluation of matrix-matrix products ${\bm S} {\bm H}$. Nevertheless, extra savings are made since ${\bm L}$ and ${\bm \Gamma}_\text{d}^{-1/2}$ acting upon matrices with $p+l$ columns can be evaluated more efficiently than those matrices ${\bm L}$ and ${\bm \Gamma}_\text{d}^{-1/2}$ multiplied with more than $p$ vectors individually in sequence when using Lanczos iteration. 

Evaluating ${\bm S}_\text{D}^T {\bm H} = {\bm L}^T \left[ {\bm S}^T \left( {\bm \Gamma}_\text{d}^{-T/2} {\bm H} \right) \right]$, where ${\bm H}$ is a matrix with $p+l$ columns, proceeds similarly to finding ${\bm S}_\text{D}^T$ times a vector with ${\bm S}^T$ times ${\bm \Gamma}_\text{d}^{-T/2} {\bm H}$ evaluated using the adjoint method (see Appendix B).

As shown in Appendices A and B, ${\bm S}$ times an $N_m \times (p+l)$ matrix and ${\bm S}^T$ times an $N_d \times (p+l)$ matrix lead to linear problems with $p+l$ right-hand sides at every time level of the direct or adjoint methods. We can take advantage of this structure and use linear solvers designed for efficient solution of linear problems with multiple right-hand sides. 

Computational time can be reduced by applying parallel solvers to solve \eqref{eq:dirnatstate}, \eqref{eq:dirprod1}, \eqref{eq:dirprodk} and (\ref{eq:adjlasttime}--\ref{eq:adjnatstate}). Another option for smaller problems, where an LU factorization is possible, is to reduce the cost of solving each right-hand side by applying direct linear solvers and re-using LU factors.

Further benefits are that the information contained in the matrices ${\bm A}^k$, ${\bm G}^k$ and ${\bm C}^k$ (see Appendix A) is only required once for each time level when finding ${\bm S} {\bm H}$ with the direct method and likewise when finding ${\bm S}^T {\bm H}$ with the adjoint method. This is unlike the Lanczos approach which requires regenerating ${\bm A}^k$, ${\bm G}^k$ and ${\bm C}^k$ or their actions on vectors over $2p$ times for every time level.

For the problems we have tested, steps \ref{algline:2viewYmat} and \ref{algline:2viewBmat} in Algorithm \ref{alg:2view} are by far the most computationally expensive since they are related to the expensive reservoir simulation. The cost of generating the random Gaussian matrix, the orthonormalization and taking the full SVD of the matrix ${\bm B}^T$ was negligible in comparison. However, for very large problems taking the full SVD of ${\bm B}^T$ can become costly. \citet{voronin2015} suggested improvements to the basic 2-view method for this situation.

Randomized methods work especially well for low-rank approximation of matrices which have rapidly decaying singular spectra. Sensitivity matrices of many inverse problems correspond to this situation since the problems are commonly ill-posed \citep{buithanh2012,vogel1993}. For cases where the singular spectrum decays slowly the accuracy of randomized methods can be improved by increasing the oversampling parameter $l$. If that fails for reasonably small $l$, then power iteration can for example be applied to improve the randomized approximation of the range of ${\bm S}_\text{D}$ or ${\bm S}_\text{D}^T$ \citep{halko2011structure,martinson2016rsvd}. However, the application of power iteration necessitates accessing ${\bm S}_\text{D}$ more often (twice for every power iteration), which negates some of the possible computational savings that can be made by applying randomized methods. Power iteration or other ways of making randomized methods more robust by increasing the number of matrix views are not considered here.

\subsubsection{Randomized 1-View Method} \label{sec:1view}
When applying the 2-view method, the adjoint simulations cannot be run at the same time as the direct simulations. This is because the randomized 2-view method requires accessing the matrix of interest twice. Randomized 1-view algorithms have been developed which aim to speed up low-rank approximation of large matrices by only using one matrix access \citep{halko2011structure,martinson2016rsvd,tropp2016,woolfe2008}. 

Here we use the 1-view approach proposed by \citet{tropp2016}. Algorithm \ref{alg:1view} presents a 1-view randomized method for estimating a TSVD based on the work of \citet{tropp2016}. Lines \ref{algline:1viewRandMats} to \ref{algline:1viewYandZ} are from the randomized sketch Algorithm 1 in \citet{tropp2016}. Lines \ref{algline:1viewQR1} to \ref{algline:1viewXmat} are from Algorithm 4 in \citet{tropp2016}. The last steps are taken from Algorithm 5 in \citet{tropp2016}. The building blocks of the 1-view method are very similar to those of the 2-view method.

\begin{algorithm}
  \caption{Randomized 1-view method for TSVD}\label{alg:1view}
  \textbf{INPUT:} {Matrix ${\bm A} \in {\mathbb{R}}^{n_r \times n_c}$ ($n_r \geq n_c$), integers $p > 0$ and $l_2 \geq l_1 \geq 0$}.\\
  \textbf{RETURNS:} {Approximate rank-$p$ SVD, ${\bm U}_p {\bm \Lambda}_p {\bm V}_p^T$, of ${\bm A}$.}
  \begin{algorithmic}[1]	
    \State Generate Gaussian random matrices ${\bm \Omega} \in {\mathbb{R}}^{n_c \times (p+l_1)}$ and ${\bm \Psi} \in {\mathbb{R}}^{(p+l_2) \times n_r}$.   \label{algline:1viewRandMats}
    \State Optional orthogonalization: ${\bm \Omega}=\text{orth}({\bm \Omega})$ and ${\bm \Psi}^T=\text{orth}({\bm \Psi}^T)$.
    \State Form the matrices ${\bm Y} = {\bm A} {\bm \Omega}$   and   ${\bm Z}^T = {\bm A}^T {\bm \Psi}^T$.    \Comment{${\bm Y} \in {\mathbb{R}}^{n_r \times (p+l_1)}$ , ${\bm Z} \in {\mathbb{R}}^{(p+l_2) \times n_c}$}    \label{algline:1viewYandZ}
    \State Find an orthonormal matrix ${\bm Q}\in {\mathbb{R}}^{n_r \times (p+l_1)}$, using QR factorization, such that ${\bm Y}={\bm Q} \tilde{\bm R}$. \label{algline:1viewQR1}
    \State Using QR factorization find matrices ${\hat{\bm Q}} \in {\mathbb{R}}^{(p+l_2) \times (p+l_1)}$ and ${\hat{\bm R}} \in {\mathbb{R}}^{(p+l_1) \times (p+l_1)}$ such that ${\bm \Psi} {\bm Q} = {\hat{\bm Q}} {\hat{\bm R}}$. 
    \State Find ${\bm X}\in {\mathbb{R}}^{(p+l_1) \times n_c}$ such that ${\hat{\bm R}} {\bm X} = {\hat{\bm Q}}^T {\bm Z}$ or ${\bm X} = {\hat{\bm R}}^{-1} {\hat{\bm Q}}^T {\bm Z}$.    \label{algline:1viewXmat}
    \State Calculate the SVD of the relatively small matrix ${\bm X}=\hat{\bm U}_{p+l_1} {\bm \Lambda}_{p+l_1} {\bm V}_{p+l_1}^T$ and truncate.
    \State Form the matrix ${\bm U}_p = {\bm Q} \hat{\bm U}_p$.
  \end{algorithmic}
\end{algorithm}

Applying the 1-view Algorithm \ref{alg:1view} to estimating a rank-$p$ SVD of ${\bm S}_\text{D}$ requires the same number of adjoint and direct solves as the 2-view method if the oversampling parameters are chosen as $l_2 = l_1 = l$ (not advised). However, unlike the 2-view approach, the adjoint solves can be run at the same time as the direct ones, see line \ref{algline:1viewYandZ} in Algorithm \ref{alg:1view} and Fig. \ref{fig:1}. 

By using a 1-view method, assuming no communication overhead and ideal parallelizability we can in theory attain similar computational speed per iteration as the BFGS method, which uses one adjoint solve to evaluate the gradient, while at the same time maintaining a convergence rate similar to standard LM. The drawback of a faster SVD approximation afforded by the 1-view method is however that it is less accurate than the 2-view approach due to the additional approximations that are made to allow for a single access approach [Tropp et al., 2016].

Note that choosing $l_2 = l_1$ is not advised as Algorithm \ref{alg:1view} is especially fragile for this choice \citep{tropp2016}. Here we chose $l_1=l=10$, for a simpler comparison with the 2-view method, and $l_2=20$, for robustness. However, for future implementations we may consider varying $l_1$ and $l_2$ adaptively based on suggestions and theoretical considerations given in \citet{tropp2016}.

\subsubsection{Subspace Iteration Method of Vogel and Wade}
\citet{vogel1994} presented a block-based subspace iteration method, which has very similar attributes to the randomized 2-view method, for evaluating a TSVD of a matrix (see Algorithm \ref{alg:subiter} in Appendix C). They presented the subspace iteration method as an alternative to the Lanczos approach for TSVD-based inversion.

The subspace algorithm proposed by \citet{vogel1994} is an adaptation of subspace methods used for eigendecomposition. Their subspace algorithm, like the randomized TSVD methods, involves computing the action of an input matrix and its transpose with thin matrices. \citet{vogel1994} recognized that this block characteristic of the subspace iteration method can be advantageous for truncated inversion using parallel architectures. Their subspace algorithm includes iteration to improve TSVD estimates. 

Appendix C elaborates on the similarities between the subspace iteration and the randomized methods. There it is shown that a modern randomized TSVD method can be obtained by making minor changes to the subspace iteration method of \citet{vogel1994}.

Like the randomized TSVD methods presented above, \citet{vogel1994} used randomized matrices to initialize the subspace iteration method for their test problems. Randomized matrices have good theoretical properties and work well in practice \citep{halko2011structure,tropp2016}. But some other initialization matrix can also be chosen if there is reason to believe that the subspace spanned by its column vectors aligns well with the principal right-singular vectors of the target matrix whose TSVD is sought.

The initialization matrix can for instance be generated by re-using the subspace spanned by the approximate right-singular vectors of a matrix closely related to the target one \citep{vogel1994}. \citet{vogel1993} and \citet{vogel1994} suggested that subspace re-use may be applicable when seeking TSVD approximations of the sequence of sensitivity matrices which arise during LM iteration.

\citet{vogel1993} found that subspace re-use can be beneficial for inversion using a truncated LM method. \citet{vogel1993} reported finding LM updates based on a truncated eigendecomposition of the GN normal matrix found using subspace iteration. They reported that one subspace iteration was often enough when using the approximate eigendecomposition from the previous LM iteration to initialize their subspace method for the following LM iteration. 

The subspace iteration method as presented by \citet{vogel1994} requires at least one more matrix access than the 2-view Algorithm \ref{alg:2view}. The subspace iteration method is therefore not applied here.
However, motivated by the subspace ideas presented in \citet{vogel1993} and \citet{vogel1994} the next section looks at re-using the subspaces spanned by the right- and left-singular vectors estimated at previous LM iterations to improve the TSVD-LM approach when applying the 1- and 2-view approaches.

\subsubsection{Subspace Re-use} \label{sec:subreuse}
Subspace ideas can be used to improve low-rank matrix approximations using randomized methods. \citet{gu2015}, for instance, proposed an improved randomized method for low-rank SVD approximation which uses the subspace framework (see Algorithm 8.1 in \citet{gu2015}). The first stage of his algorithm finds a low-rank SVD of the input matrix using a randomized 2-view method. The second stage of the algorithm uses the subspace spanned by the approximate right-singular vectors, found at the first stage, to initialize a power or subspace iteration method which outputs an improved low-rank approximation. 

Another simpler way of improving the TSVD approximation is to use the 2-view method again at the second stage instead of the power/subspace iteration. This type of subspace re-use can likewise improve low-rank matrix approximations when using the 1-view method.

This suggests that a subspace re-use scheme like that used by \citet{vogel1993} might also be reasonable when running the TSVD-LM method with the 1- or 2-view approaches. That is using the singular vectors from a previous LM iteration to initialize the 2-view or 1-view methods at the following LM iteration.

Re-using the subspace spanned by the singular vectors from the previous LM iteration may work well if the sensitivity matrix does not change much between iterations. This may be the case for problems that are close to linear or at late LM iterations where the model parameter updates are small. 

The results presented by \citet{tavakoli2010,shirangi2011,shirangi2014} indicate that the singular spectrum of ${\bm S}_\text{D}$ is reasonably stable during their inversion runs, especially at late LM iterations. We therefore speculate that the randomized methods with subspace re-use may work well for inversion of some reservoir problems. This approach may, however, prove to be ineffective for geothermal problems which are usually very nonlinear. Furthermore, subspaces from a previous LM iteration may be a bad choice for "large" model updates. Nevertheless, use of appropriate randomized oversampling may prevent failure of the subspace approach.

Based on these ideas, we developed versions of Algorithms \ref{alg:2view} and \ref{alg:1view} that apply subspace re-use. These variants are considered in an attempt to improve the basic randomized methods. The cost of the re-use variants remains nearly the same as the original algorithms when used within the TSVD-LM scheme. The only modifications have to do with the sampling matrices ${\bm \Omega}$ and ${\bm \Psi}$. 

Running the TSVD-LM method using the randomized 2-view method with subspace re-use applies Algorithm \ref{alg:2view} as is at the first LM iteration. For subsequent LM iterations the first $p^\text{prev}$ columns of ${\bm \Omega}$ are given by the $p^\text{prev}$ approximate right-singular vectors of ${\bm S}_\text{D}$ retained at the previous iteration, when the 2-view method is applied to ${\bm S}_\text{D}$. That is
\begin{linenomath*}
\begin{equation}
	{\bm \Omega} = 
    \begin{bmatrix}
    {\bm v}_1 & \dots & {\bm v}_{p^\text{prev}} & {\bm \omega}_1 & \dots & {\bm \omega}_{l+p-p^\text{prev}}
    \end{bmatrix}
\, ,
\end{equation} \end{linenomath*}
where ${\bm \omega}_i$ denotes a random column vector. However, when applying the 2-view method to ${\bm S}_\text{D}^T$, the first $p^\text{prev}$ columns of ${\bm \Omega}$ are determined by the retained left-singular vectors of ${\bm S}_\text{D}$ from the previous iteration. Subspace re-use may therefore reduce the number of random vectors ${\bm \omega}_i$ generated at a LM iteration, which is a possible benefit.

The randomized 1-view method using subspace re-use treats ${\bm \Omega}$ in the same way as outlined for the 2-view method with subspace re-use and ${\bm \Psi}^T$ is treated in the same way as the case of ${\bm \Omega}$ applied to the input matrix ${\bm A}$ transposed.

The following sections compare the performance of the TSVD-LM method using the standard Lanczos approach with TSVD-LM run using the 1-view and 2-view methods, with and without subspace re-use.

\section{Description of Computational Experiments} \label{sec:expdesign}
\subsection{Comparison of TSVD Methods}
The following computational experiments were carried out to demonstrate the applicability of the proposed randomized TSVD-LM methods. The study looks at comparing: (i) methods used to control the number of retained singular values $p$ (increasing $p$ linearly or varying sv-cut), (ii) the computational efficiency and convergence characteristics of inversions using randomized methods against Lanczos based inversions, (iii) the robustness of the 1-view and 2-view methods, and (iv) randomized methods with and without subspace re-use.

Item (i) is of interest since the randomized implementations increase the truncation level $p$ linearly between LM iterations unlike the Lanczos methods tested by \citet{shirangi2014,shirangi2016,tavakoli2010,tavakoli2011}, which used the sv-cut approach. The suitability of the linear truncation approach is investigated in section \ref{sec:resultstrunccontrol} by comparing TSVD-LM inversions using the Lanczos method with linear and sv-cut controlled truncation. Items (ii)-(iv) are addressed in section \ref{sec:resultsLancVsRand} by comparing inversion results applying the standard Lanczos approach against the four randomized methods under consideration. Table \ref{tab:methods} lists the attributes of the six methods compared in this study. 

\begin{table}[h!]
\centering
\caption{The six TSVD approaches applied in the study.}
\label{tab:methods}      
\begin{tabular}{lllll}
\hline\noalign{\smallskip}
\textbf{Method} & \textbf{Subspace} & \textbf{Algorithm} & \textbf{Truncation} &  \textbf{Oversampling}  \\
 & \textbf{Re-use} &  & \textbf{Control} &  \textbf{Parameter}  \\
\noalign{\smallskip}\hline\noalign{\smallskip}
Lanczos & No & \ref{alg:lanczos} & Linear & $\varepsilon_\text{sv} = 10^{-5}$ or $10^{-1}$ \\
        & No & \ref{alg:lanczos} & sv-cut & $\varepsilon_\text{sv} = 10^{-5}$ \\
2-view & No  & \ref{alg:2view} & Linear & $l = 10$ \\
       & Yes & \ref{alg:2view}$^\textnormal{a}$ & Linear & $l = 10$ \\
1-view & No  & \ref{alg:1view} & Linear & $l_1 = 10$ and $l_2 = 20$ \\
       & Yes & \ref{alg:1view}$^\textnormal{a}$ & Linear & $l_1 = 10$ and $l_2 = 20$ \\
\noalign{\smallskip}\hline
\multicolumn{5}{l}{ $^\textnormal{a}$With subspace modifications, see end of section \ref{sec:subreuse}.}
\end{tabular}
\end{table}

Computational experiments were performed using a two-dimensional synthetic high enthalpy pure water reservoir model, previously constructed to test the applicability of adjoint and direct methods to inversion of geothermal models \citep{bjarkason2016}. Reservoir simulations were run using AUTOUGH2 \citep{angusyeh2012}, The University of Auckland's version of TOUGH2 \citep{pruess2004,pruess1999tough2}. The methods described in this work were implemented using Python scripts and experiments were run on a 3.40 GHz Intel i7-4770 CPU with 16 GB RAM and 4 cores. The Numba JIT compiler was used to speed up parts of the Python code \citep{lam2015numba}.

\subsection{Synthetic Truth Model}
\begin{figure}
\centering
\includegraphics[width=0.95\textwidth]{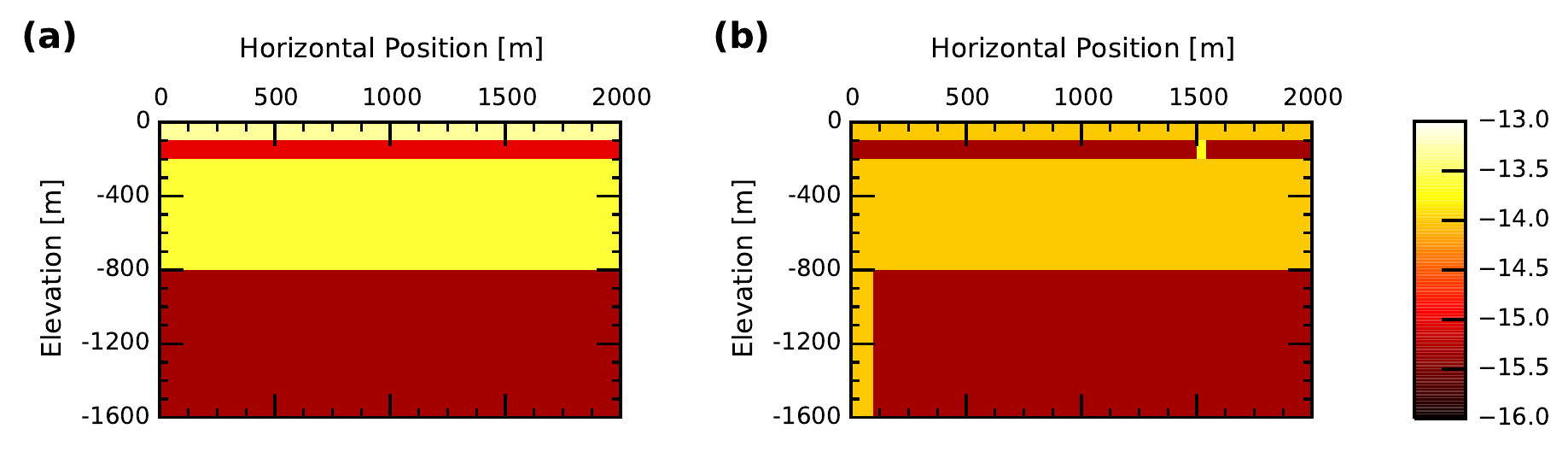}
\caption{True formation permeabilities: (a) logarithm of the horizontal permeability $k_x$ [m$^2$]; (b) logarithm of the vertical permeability $k_z$ [m$^2$].}
\label{fig:2}  
\end{figure}

The two-dimensional synthetic vertical slice model used here covers a vertical $1.6$~km by $2$~km rectangular area and consists of $8,100$~blocks, see Fig. \ref{fig:2}. The main body of the reservoir is made up of $8,000$ equally sized $20$~m$\times 20$~m$\times 20$~m blocks. The additional $100$ blocks are large volume blocks at the top boundary that effectively keep it at a constant $15^\circ$C temperature and a pressure of $101.35$~kPa, while the lateral side boundaries are closed. Convection is induced in the system by a constant mass flux of $0.15$~kg/s, with an enthalpy of $1,200$~kJ/kg, spread evenly over the first $100$~m of the bottom boundary. A constant heat flux of $80$~mW/m$^2$ is applied over the remainder of the bottom boundary.

The true formation permeabilities are shown in Fig. \ref{fig:2}. Here the task is to infer these permeabilities using synthetic data.
More details can be found on the model in \citep{bjarkason2016}.

\subsection{Synthetic Data}
\begin{figure}
\centering
\includegraphics[width=0.95\textwidth]{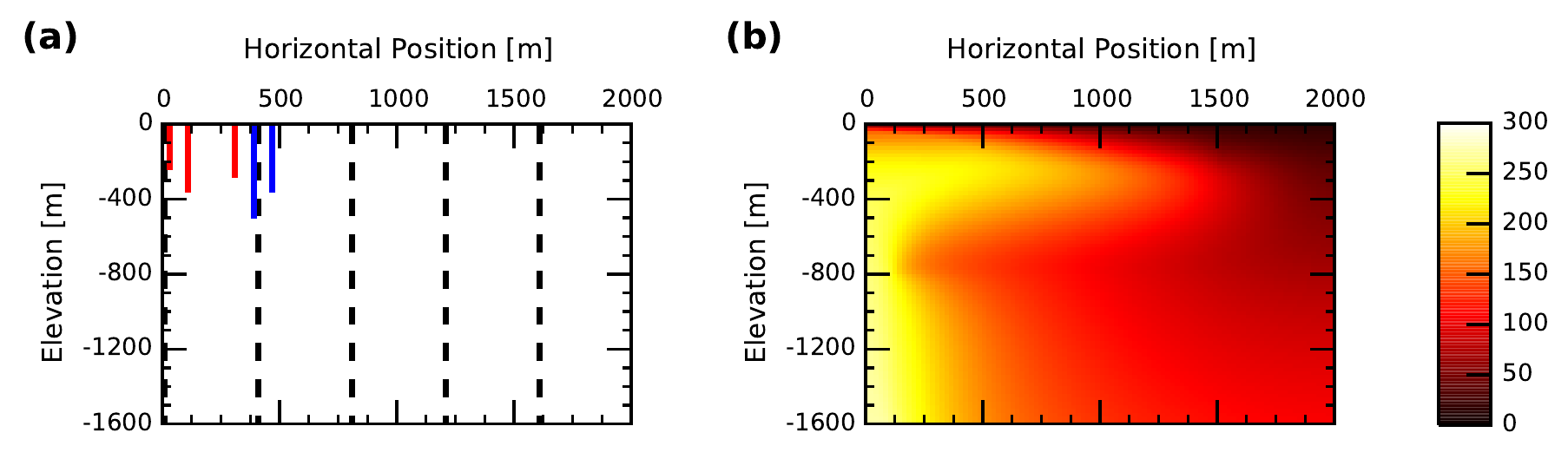}
\caption{(a) Location of the 3 production wells (red), 2 injection wells (blue) and 5 natural state observation wells (black dashed lines). (b) True natural state temperatures [$^\circ$C].}
\label{fig:3}  
\end{figure}

Synthetic observations were generated by running a natural state simulation for the true model followed by a production run. The natural state simulations use a transient approach to find the steady natural state conditions that exist prior to the production period. The true natural state temperature distribution is shown in \ref{fig:3}(b). Natural state temperatures were recorded for every third model block intersected by the black dashed well tracks shown in Fig. \ref{fig:3}(a), giving 135 temperature observations. Gaussian random noise with a standard deviation of $0.5^\circ$C was added for the final temperature observations.

The production scenario lasts three years. During production, mass is extracted at constant rates from three production wells and some of the fluid is reinjected back into the system at two injection wells. Figure \ref{fig:3}(a) shows the well tracks of the production and injection wells. Numbering the production wells in Fig. \ref{fig:3}(a) from left to right, both Producer 1 and 2 extract $0.5$~kg/s, while Producer 3 extracts $0.4$~kg/s. Reinjection of $1.12$~kg/s of $167.5$~kJ/kg enthalpy water is divided equally between the two injection wells. The production and injection feed zones coincide with the bottom of the respective well tracks shown in Fig. \ref{fig:3}(a).

Production history observations are production pressures and enthalpies at the feed zones of the three production wells. Taking measurements every three months at all production wells gave 36 pressure observations and 36 enthalpy observations. Gaussian noise was also added to the production data, using standard deviations of $0.1$~bar and $10$~kJ/kg for pressures and enthalpies, respectively.

This gave a total of $N_d = 207$ observations, which is rather small compared to the number expected for a real geothermal reservoir, but suffices for comparing the randomized methods against the Lanczos approach.

\subsection{Inversion Setup}
The synthetic observations were used for inverting the log-transformed horizontal and vertical permeabilities of every model block, apart from the top boundary blocks, using variants of the TSVD approaches, see Table \ref{tab:methods}. Inversions, therefore, included $16,000$ parameters. The vector of adjustable parameters was ${\bm m} = [\log{k_{x,1}} \, , \dots \, , \log{k_{x,N_m/2}} \, ,  \log{k_{z,1}} \, , \dots \,$ $ , \log{k_{z,N_m/2}}] $, where $\log{k_{z,i}}$ and $\log{k_{x,i}}$ are the base-ten logarithms of the vertical and horizontal permeabilities, respectively, for the $i$th adjustable model block. Permeabilities were allowed to vary between $10^{-16}$ and $10^{-13}$~m$^2$. All inversions began with all adjustable parameters set equal to ${\bm m}_\text{pr}=-14$.

The inversions used a regularization scheme similar to that used by \citet{shirangi2011} for investigating the use of the TSVD-LM method of \citet{tavakoli2010} for non-Gaussian parameter fields. The regularization scheme used here was designed to give preference to models with spatially smooth permeability valeus, and locally similar horizontal and vertical permeabilities. The regularization matrix is chosen here as ${\bm R} = {\bm W}^T {\bm W}$ where
\begin{linenomath*}
\begin{equation} \label{eq:Wmat}
	{\bm W} = 
    \begin{bmatrix}
    	{\bm L}_{k_x} \\
        {\bm L}_{k_z} \\
        {\bm L}_{{k_x}{kz}} \\
        10^{-3} {\bm I}
    \end{bmatrix}
    \, .
\end{equation} \end{linenomath*}
Here ${\bm L}_{k_x} = [ {\bm L}_1 \, {\bm 0} ]$ and ${\bm L}_{k_z} = [ {\bm 0} \, {\bm L}_1 ]$ induce spatial smoothing, where ${\bm L}_1 \in {\mathbb{R}}^{N_\text{conadj} \times N_m/2}$ is a discrete representation of the first derivative operator and $N_\text{conadj}$ is the number of connections between adjustable model blocks. For connection number $c$ connecting blocks $i$ and $j$ we have $[{\bm L}_1 ]_{ci} = 1$ and $[{\bm L}_1 ]_{cj} = -1$, all other elements of ${\bm L}_1$ are zero. ${\bm L}_{k_x k_z } \in {\mathbb{R}}^{N_m/2 \times N_m }$ is a sparse matrix with $[{\bm L}_{k_x k_z } ]_{ii}=1$ and $[{\bm L}_{k_x k_z } ]_{i(i+N_m/2)} = -1$ which suggests that the horizontal and vertical permeabilities should be similar. The identity term in \eqref{eq:Wmat} is included to ensure that ${\bm R}$ is positive definite so we can apply Cholesky factorization. 

The matrix ${\bm L}$ is defined such that ${\bm R}^{-1} = {\bm L} {\bm L}^T$, however, to avoid inverting ${\bm R}$ we find the Cholesky factorization of ${\bm R} = {\bm L}^{-T} {\bm L}^{-1}$ and store ${\bm L}^{-1}$ instead of ${\bm L}$.

Assuming no modeling errors apart from those introduced by the observation noise, then sensible models can be expected to result in observation mismatch terms with \citep{oliver2008,tarantola2005}
\begin{linenomath*}
\begin{equation} \label{eq:obsmismespect}
	105 = N_d - 5 \sqrt{2 N_d} \leq \Phi_\text{d} ( \bm{m} ) \leq  N_d + 5 \sqrt{2 N_d} = 309
\end{equation} \end{linenomath*}
with the expectation value of $\Phi_\text{d} ( \bm{m} ) = N_d$. Introducing the normalized observation mismatch $\Phi_\text{N} ( \bm{m} ) = \Phi_\text{d} ( \bm{m} ) / N_d$, \eqref{eq:obsmismespect} can be written as
\begin{linenomath*}
\begin{equation}
	0.5 = 1 - 5 \sqrt{2/N_d} \leq \Phi_\text{N} ( \bm{m} ) \leq  1 + 5 \sqrt{2/N_d} = 1.5
\end{equation} \end{linenomath*}
Experiments showed that a regularization weight $\mu = 2.5$ usually resulted in $\Phi_\text{d} ( \bm{m} ) \approx N_d$ for the present problem. To simplify the comparison of the various TSVD-LM methods all results presented here used a fixed regularization weight $\mu = 2.5$.

\section{Results of Computational Experiments} \label{sec:results}
The TSVD-LM method was used to solve the inverse problem described in the previous section. Inversions were run using the Lanczos method as well as the randomized 1- and 2-view methods. All methods used the direct linear solver SuperLU to solve the adjoint and direct equations to estimate the TSVD of ${\bm S}_\text{D}$. To compare the methods, inversions were halted after 30 LM iterations, unless stated otherwise.

\subsection{Singular Spectrum of Dimensionless Sensitivity Matrix}
\begin{figure}
\centering
\includegraphics[width=0.99\textwidth]{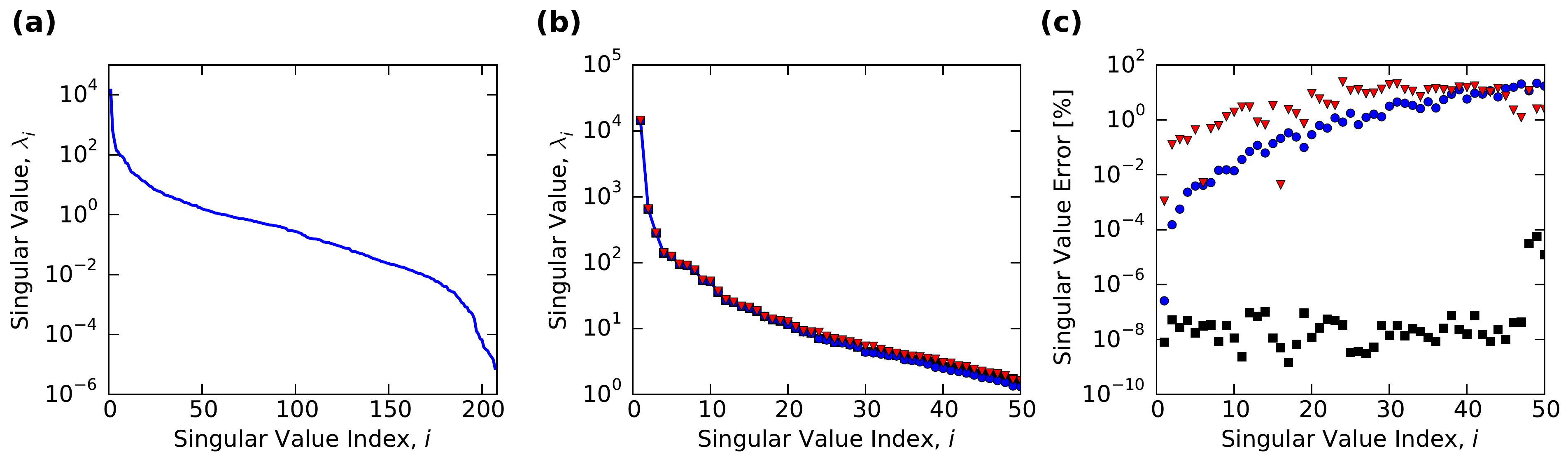}
\caption{For initial parameters and truncation $p=50$: (a) Full singular spectrum, (b) $50$ largest singular values, (c) estimation errors. (Solid blue line) Actual singular values, (black squares) values using Lanczos iteration ($\varepsilon_\text{sv} = 10^{-5}$), (blue circles) values using 2-view method and (red triangles) values using 1-view method.}
\label{fig:4}     
\end{figure}

Figure \ref{fig:4} shows for the initial parameter guess the actual singular spectrum of ${\bm S}_\text{D}$ (found by forming the full matrix ${\bm S}_\text{D}$ using the adjoint method) and estimates of the $50$ largest singular values using the Lanczos ($\varepsilon_\text{sv} = 10^{-5}$), the 1-view ($l_1=10$ and $l_2=20$), and the 2-view ($l=10$) methods. As expected the results show that the 1-view method is less accurate than the 2-view method. The singular spectrum shows an initial rapid decay of singular values but flattens out for smaller singular indices. The relative flatness of the spectrum results in inaccurate estimates of the small singular values using the randomized methods. The Lanczos method is a lot more accurate than the randomized methods.

\subsection{Truncation Control} \label{sec:resultstrunccontrol}
\subsubsection{Linearly Adjusted Truncation Compared with sv-cut Approach} 
Inversion was initially performed using the Lanczos method, with a convergence tolerance $\varepsilon_\text{sv} =10^{-5}$. The computational cost was reduced by gradually increasing the number of retained singular values $p$ between successive LM iterations, increasing $p$ linearly or controlling $p$ using sv-cut. 

The linear approach used $p=1$ at the first LM iteration and increased $p$ by two between LM iterations. Like \citet{shirangi2014} and \citet{shirangi2016} the sv-cut approach used $\text{sv-cut}=0.5$ for the first LM iteration and sv-cut was halved between LM iterations. In all cases $p$ was not allowed to exceed $50$, which is a similar maximum value to that used by \citet{shirangi2014} and \citet{shirangi2016}.

\begin{figure}
\centering
\includegraphics[width=0.95\textwidth]{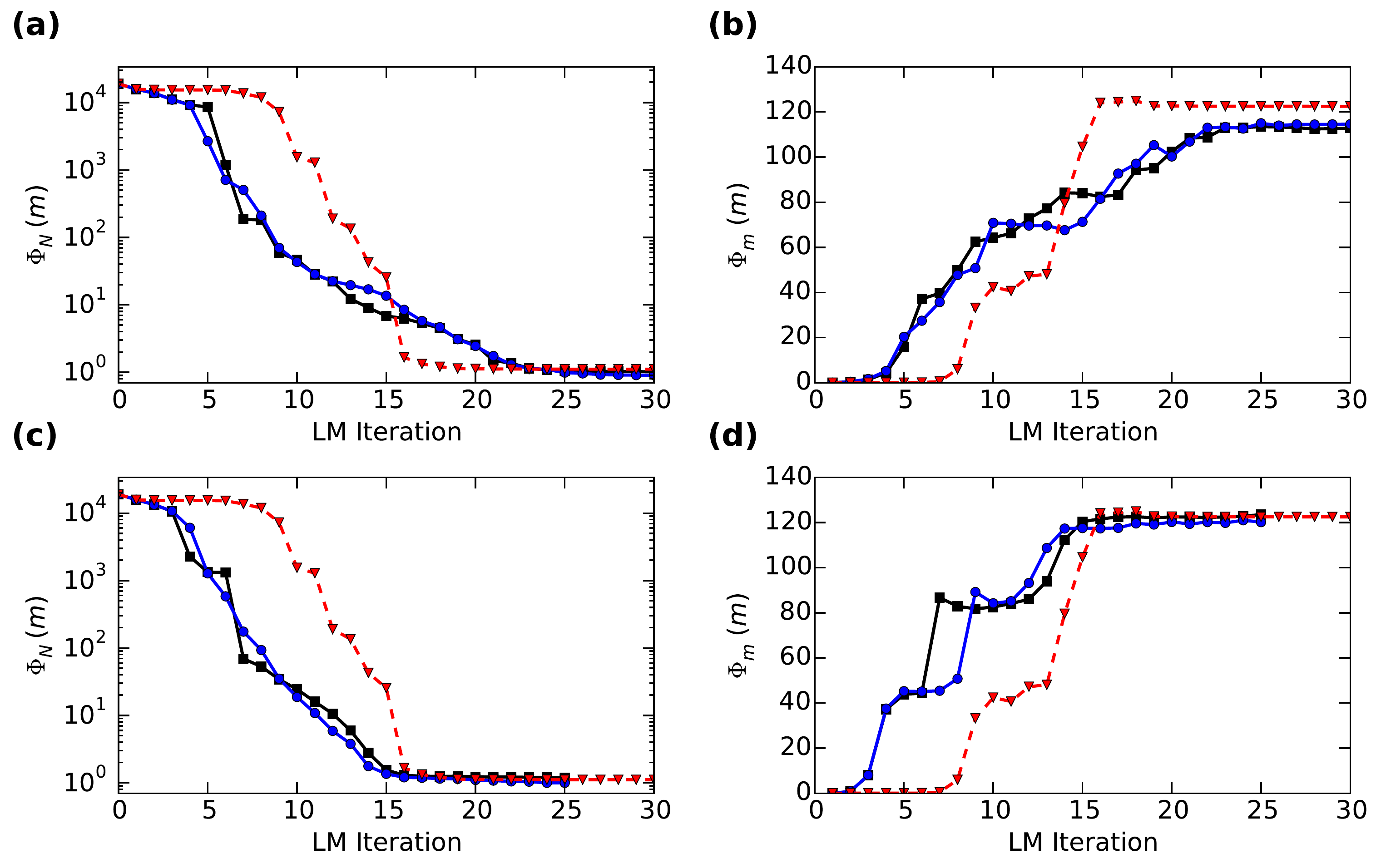}
\caption{(a) Normalized observation mismatch and (b) model mismatch using TSVD-LM with Lanczos iteration. (Black squares) Truncation $p$ increased linearly by two between LM iterations ($\varepsilon_\text{sv} =10^{-5}$), (blue circles) truncation $p$ increased linearly by two between LM iterations ($\varepsilon_\text{sv} =10^{-1}$) and (red triangles) truncation adjusted using sv-cut ($\varepsilon_\text{sv} =10^{-5}$). (c)-(d) Show the same as (a)-(b) but with $p$ increased by five between LM iterations for the linear increase approach.}
\label{fig:5}       
\end{figure}

Figures \ref{fig:5}(a) and \ref{fig:5}(b) compare the convergence performance of the linear truncation and sv-cut approaches in terms of the normalized observation mismatch and the model mismatch. The linear approach resulted in a final observation mismatch $\Phi_\text{d} = 210$ and the sv-cut approach resulted in a slightly worse match with $\Phi_\text{d} = 230$. Both matches comply with the bounds given in \eqref{eq:obsmismespect}. The results in Fig. \ref{fig:5}(b) show that the sv-cut approach resulted in a rougher model with a larger model mismatch term. 

\begin{figure}
\centering
\includegraphics[width=0.9\textwidth]{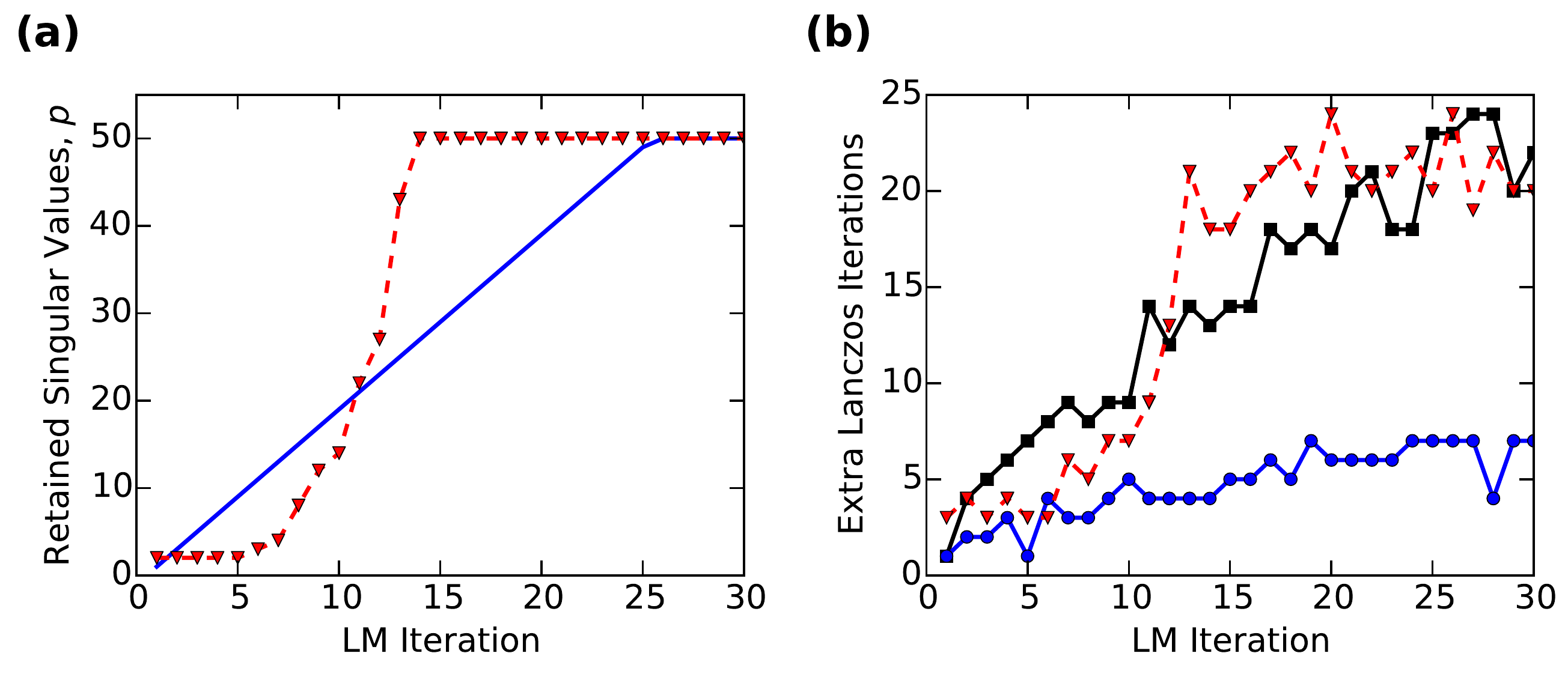}
\caption{ (a) (Solid blue line) Number of retained singular values $p$ increased linearly by two between LM iterations and (red triangles) number of singular values when using sv-cut. (b) Number of Lanczos iterations in excess of the truncation $p$. (Black squares) $p$ increased linearly ($\varepsilon_\text{sv} =10^{-5}$), (blue circles) $p$ increased linearly ($\varepsilon_\text{sv} =10^{-1}$) and (red triangles) $p$ adjusted with sv-cut ($\varepsilon_\text{sv} =10^{-5}$).}
\label{fig:6}     
\end{figure}

Figure \ref{fig:6}(a) depicts the number of retained singular values $p$ for the two Lanczos inversion approaches. The sv-cut approach resulted in a near step-like increase of $p$. The initial slow increase of $p$ resulted in an initially slow decrease of the observation mismatch, see Fig. \ref{fig:5}(a). The rapid increase in $p$ midway through the inversion resulted in a sudden rapid decrease in the objective function. The rapid increase in $p$ when applying the sv-cut method is because of the gradual flattening of the singular spectrum. This rapid increase in the number of retained singular values may have contributed to the increased model roughness found using the sv-cut approach.

The objective function can be lowered more rapidly when using the linear truncation control method by increasing $p$ by a larger value between LM iterations. This can be seen in the results shown in Figs. \ref{fig:5}(c) and \ref{fig:5}(d), obtained by increasing $p$ by five between iterations. The inversion was limited to $25$ LM iterations. Again, the rapid increase in the number of singular values resulted in a rougher model. 

When experimenting with the TSVD-LM method, we generally found that slower increases in both the number of retained singular values $p$ and the LM damping factor $\gamma$ helps to reduce model roughness, which is consistent with the discussion in section \ref{sec:LMdamping}. 

The above results and discussion indicate that the linear truncation control approach is just as effective as the sv-cut approach. Therefore, it appears to be just as effective to use the linear truncation approach rather than sv-cut when applying the randomized TSVD-LM methods to the present problem.

\subsubsection{Reducing Cost by Increasing the Lanczos Convergence Tolerance} 
For this problem we found that using a low Lanczos convergence tolerance $\varepsilon_\text{sv}$ resulted in a large number of excess Lanczos iterations $n$. Figure \ref{fig:6}(b) shows that at initial iterations when $p$ was small the number of excess iterations were below $10$ which agrees with values reported by \citet{tavakoli2010} and \citet{shirangi2011}. However, at late iterations the excess Lanczos iterations could be more than $20$, slowing the inversion down. 

We found that a more forgiving convergence tolerance of $\varepsilon_\text{sv} =10^{-1}$ substantially reduced the number of excess iterations, see Fig. \ref{fig:6}(b), and therefore also the computational burden. The larger tolerance did not degrade the performance of the inversion process as shown in Fig. \ref{fig:5}. When using $\varepsilon_\text{sv} =10^{-1}$ the number of excess iterations were consistent throughout with those reported by \citet{tavakoli2010} and \citet{shirangi2011}.

\subsection{Randomized vs Lanczos} \label{sec:resultsLancVsRand}
\subsubsection{Convergence Comparison}
Inversions were also carried out using the randomized 2-view and 1-view methods, with and without subspace re-use. When applying the randomized methods the number of retained singular values $p$ was at first regulated in the same linear way that worked well for the Lanczos method with $p$ increased by two between LM iterations. Inversions were run $20$ times for each randomized method listed in Table \ref{tab:methods}. For all following results the randomized methods used oversampling parameters $l=l_1=10$ and $l_2=20$, and $p$ was allowed a maximum value of $50$.

\begin{figure}
\centering
\includegraphics[width=0.95\textwidth]{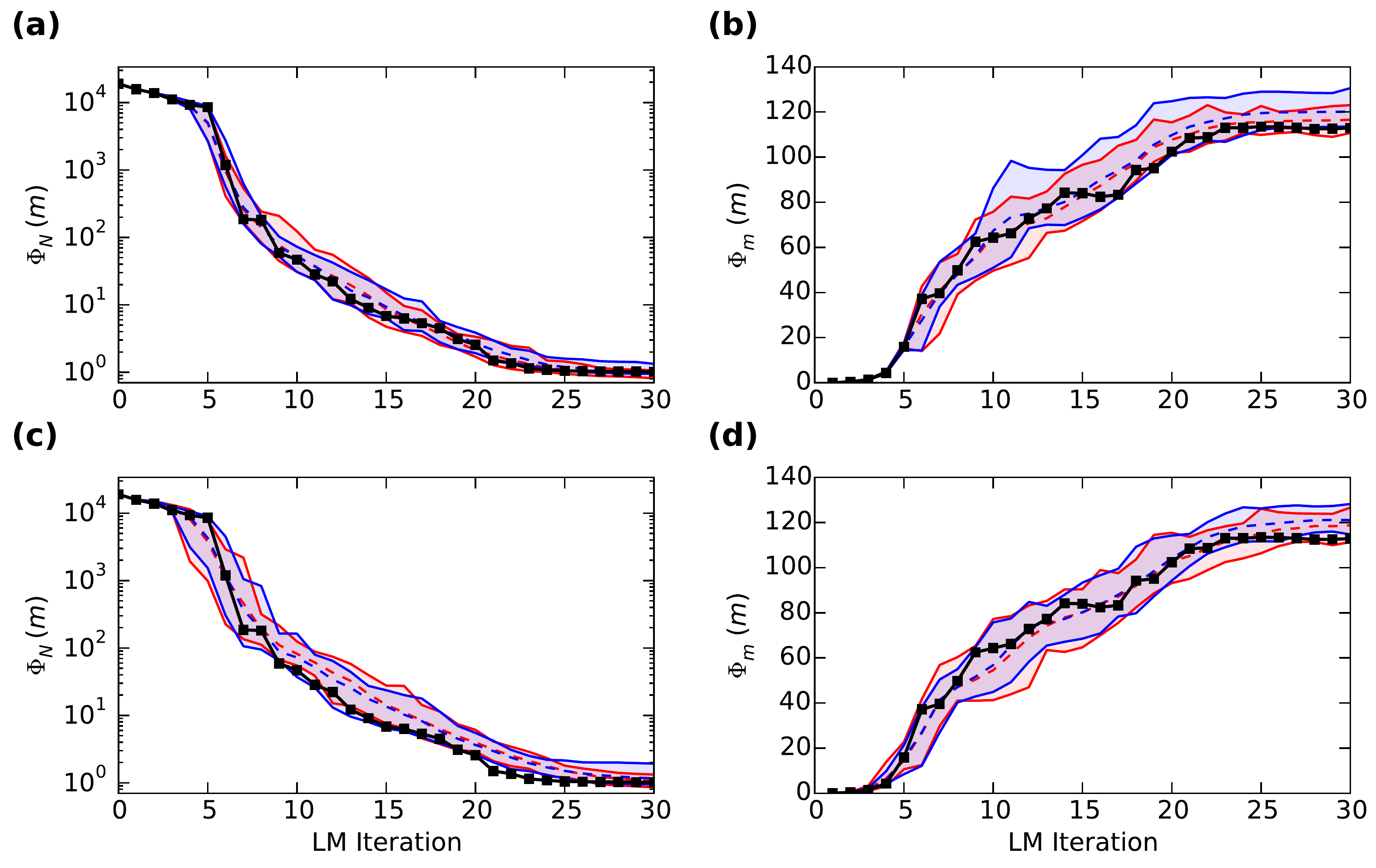}
\caption{(a) Normalized observation mismatches and (b) model mismatches using TSVD-LM with Lanczos iteration with $\varepsilon_\text{sv} =10^{-5}$ (black squares) compared with the 2-view method. The dashed red line indicates mean values when running 20 inversions with the 2-view method without subspace re-use. Red solid lines indicate the maximum and minimum values of the 20 runs. The same is plotted in blue for the 2-view method using subspace re-use. (c) and (d) show the same comparison when using the 1-view method instead of the 2-view one. For all runs the number of singular values was initially set to one and subsequently increased by two between iterations.}
\label{fig:7}   
\end{figure}

The normalized observation mismatches and model mismatches as a function of LM iteration found using the randomized methods are compared against the Lanczos method (with $\varepsilon_\text{sv} =10^{-5}$ and $p$ increased linearly by two) in Fig. \ref{fig:7}. The convergence properties of the randomized 2-view methods are very similar to the ones found using the Lanczos approach. However, the convergence of the 1-view approaches appear on average to be slightly slower, based on the LM iteration count, than the convergence of the Lanczos and 2-view approaches. Additional oversampling could improve the 1-view approach, but with the drawback of increasing the computational burden.

Table \ref{tab:results1} gives more details of the final objective function values and mismatch terms found using the randomized methods as well as the Lanczos methods with tolerances $\varepsilon_\text{sv} =10^{-5}$ and $\varepsilon_\text{sv} =10^{-1}$. The results show that the randomized methods gave similar model mismatch terms as those found using the Lanczos methods. 

The final values show that the 2-view method without subspace re-use gave very similar observation and model mismatches to the Lanczos runs, with the average 2-view run not using subspace re-use coinciding well with the values found using Lanczos iteration. The 2-view method applying subspace re-use and the 1-view methods, on the other hand, gave on average somewhat worse matches. The subspace variants generally performed worse than the methods not using subspace re-use.

It should be noted that the randomized methods became even more consistent with the Lanczos approach when the LM damping factor $\gamma$ was lowered less aggressively, an expected result since a larger $\gamma$ suppresses the influence of small singular values which are evaluated less accurately. 

\begin{table}[t!]
\centering
\caption{Results after 30 LM iterations. For the randomized methods the values shown are the mean (ave), minimum (min) and maximum (max) values found running the methods 20 times. For all runs the number of singular values was initially set to one and subsequently increased by two between LM iterations.}
\label{tab:results1} 
\begin{tabular}{llllllllll}
\hline\noalign{\smallskip}
\textbf{Method} & \multicolumn{3}{l}{\textbf{Objective}} & \multicolumn{3}{l}{\textbf{Obs. Mismatch}} & \multicolumn{3}{l}{\textbf{Model Mismatch}}   \\
 & $\Phi^\text{ave}$ & $\Phi^\text{min}$ & $\Phi^\text{max}$ & $\Phi^\text{ave}_\text{d}$ & $\Phi^\text{min}_\text{d}$ & $\Phi^\text{max}_\text{d}$ & $\Phi^\text{ave}_\text{m}$ & $\Phi^\text{min}_\text{m}$ & $\Phi^\text{max}_\text{m}$  \\
\noalign{\smallskip}\hline\noalign{\smallskip}
Lanczos$^\textnormal{a}$ & 492 &  &  & 210 &  &  & 113 &  &  \\ 
Lanczos$^\textnormal{b}$ & 473 &  &  & 187 &  &  & 115 &  &  \\ 
2-view & 477 & 453 & 503 & 185 & 168 & 220 & 117 & 111 & 123 \\ 
2-view$^\textnormal{c}$ & 526 & 480 & 602 & 226 & 196 & 276 & 120 & 114 & 131 \\ 
1-view & 518 & 456 & 580 & 221 & 178 & 274 & 119 & 111 & 127 \\ 
1-view$^\textnormal{c}$ & 544 & 509 & 705 & 241 & 196 & 398 & 121 & 115 & 128 \\ 
\noalign{\smallskip}\hline
\multicolumn{9}{l}{$^\textnormal{a}$ $\varepsilon_\text{sv} = 10^{-5}$; $^\textnormal{b}$ $\varepsilon_\text{sv} = 10^{-1}$; $^\textnormal{c}$with subspace re-use.}
\end{tabular}
\end{table}

The subspace re-use approaches gave the worst performance for the presented test case. Note that since the parameters outnumber the observations the 2-view method (Algorithm \ref{alg:2view}) uses ${\bm S}_\text{D}^T$ as the input matrix and the subspace re-use is performed in the observation space. This is unlike \citet{vogel1993} who applied the subspace re-use in the model parameter space. 

We also tried modifying Algorithm \ref{alg:2view} to instead use ${\bm S}_\text{D}$ as input for the present inverse problem. However, using the modified Algorithm \ref{alg:2view} to run the 2-view method with and without subspace re-use did not result in significantly different results to those presented here. Running the modified 2-view re-use method 20 times gave an average objective function value of about $495$, which is better than using subspace re-use in the observation space but worse than the 2-view method without subspace re-use.

\subsubsection{Time Spent on Inversions}
Table \ref{tab:results2} compares computational time needed to run inversions and the main computational expenses for the inversion runs presented in Table \ref{tab:results1} and Fig. \ref{fig:7}. The inversions using randomized methods were an order of magnitude faster than those using the Lanczos approach. Relaxing the Lanczos convergence tolerance $\varepsilon_\text{sv}$ helps to speed up the Lanczos runs, but the improvement is not close to that found using the randomized methods. 

The randomized methods are faster since the adjoint and direct linear problems at every time level (see Appendices A and B) are solved with one call to the SuperLU solver, which forms one LU factorization and re-uses the factorization for the multiple right-hand sides.

\begin{table}[h!]
\centering
\caption{Computational cost of inversions after 30 LM iterations. For the randomized methods the mean (ave), minimum (min) and maximum (max) values of the number of simulations and inversion time were found by running each method 20 times.}
\label{tab:results2}       
\begin{tabular}{lllllllll}
\hline\noalign{\smallskip}
\textbf{Method} & \multicolumn{3}{l}{\textbf{Simulations}} & \textbf{Direct} & \textbf{Adjoint} & \multicolumn{3}{l}{\textbf{Time} [s]}  \\
 & ave & min & max & & & ave & min & max \\
\noalign{\smallskip}\hline\noalign{\smallskip}
Lanczos$^\textnormal{a}$ & 36 &  &  & 1,335 & 1,305 & 68,200 &  &  \\
Lanczos$^\textnormal{b}$ & 36 &  &  & 1,047 & 1,017 & 54,500 &  &  \\
2-view & 36.5 & 33 & 41 & 1,175 & 1,175 & 7,100 & 6,750 & 7,590 \\
2-view$^\textnormal{c}$ & 35.5 & 32 & 40 & 1,175 & 1,175 & 7,000 & 6,660 & 7,410 \\
1-view & 37.3 & 33 & 41 & 1,475 & 1,175 & 6,040 & 5,620 & 6,570 \\
1-view$^\textnormal{c}$ & 35.6 & 31 & 47 & 1,475 & 1,175 & 5,870 & 5,510 & 7,380 \\
\noalign{\smallskip}\hline
\multicolumn{9}{l}{$^\textnormal{a}$ $\varepsilon_\text{sv} = 10^{-5}$; $^\textnormal{b}$ $\varepsilon_\text{sv} = 10^{-1}$; $^\textnormal{c}$with subspace re-use.}
\end{tabular}
\end{table}

Running 30 LM iterations was usually faster using the 1-view approach than the 2-view approach since the adjoint and direct solves were run in parallel when applying the 1-view approach. The 1-view method did not result in substantial speed-up since the cost of inversions running the random methods were largely dominated by the cost of the nonlinear forward simulations and the 1-view method used additional oversampling.

The total number of nonlinear simulations were fairly similar for all methods. Compared to the Lanczos method, in some cases the randomized methods resulted in fewer simulations and in some instances required an increased number of simulations. The total number of adjoint and direct runs was also similar for all methods.

\subsubsection{Estimated Parameters and Observation Matches}
\begin{figure}[b!]
\centering
\includegraphics[width=0.95\textwidth]{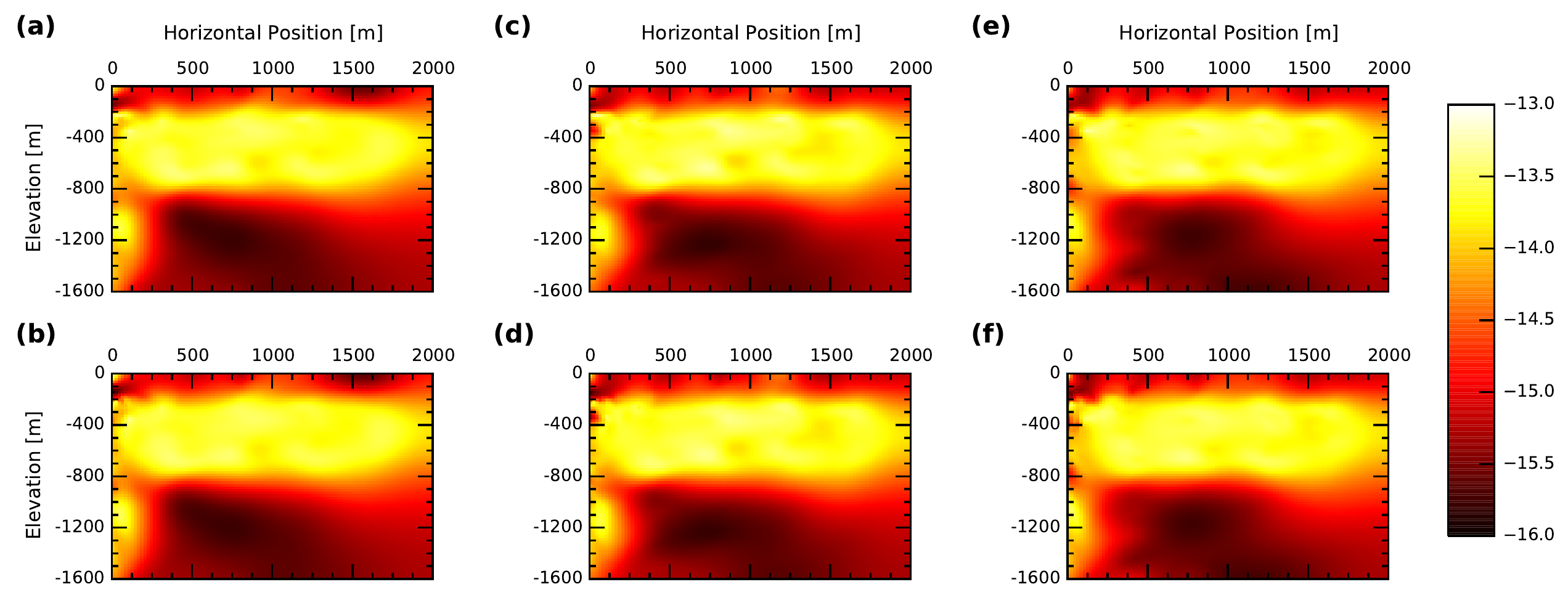}
\caption{Estimated log-transformed permeabilities, (a) horizontal and (b) vertical, using the Lanczos method ($\varepsilon_\text{sv} = 10^{-5}$). (c)-(d) and (e)-(f) show the same using the 2-view and 1-view methods, respectively, without subspace re-use.}
\label{fig:8}      
\end{figure}

Figures \ref{fig:8}(a) and \ref{fig:8}(b) show the estimated permeabilities using the Lanczos method with $\varepsilon_\text{sv} = 10^{-5}$. The inversion managed to capture the large-scale features of the true permeability field shown in Figure \ref{fig:2} as well as the permeability contrast above the feed zone of Producer 1 due to the low permeability cap-rock. Figures \ref{fig:8}(c)-(d) and \ref{fig:8}(e)-(f) show that very similar results were found for the first runs made using the randomized 2-view and 1-view methods without subspace re-use. The estimated permeability values found for other runs of the randomized methods with and without subspace re-use look similar to those shown in Figs. \ref{fig:8}(c)-(f).

\begin{figure}
\centering
\includegraphics[width=0.95\textwidth]{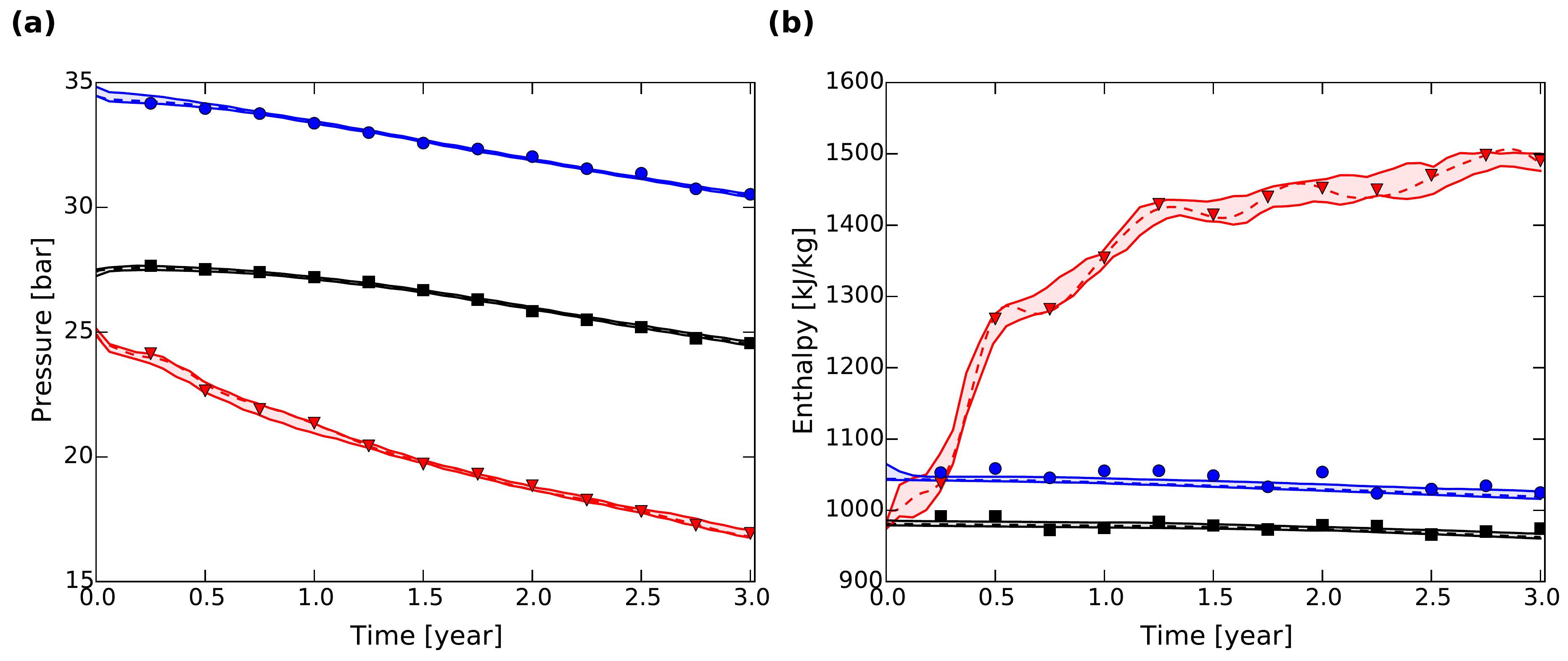}
\caption{Production (a) pressure and (b) enthalpy observations for Producer 1 (red triangles), Producer 2 (blue circles) and Producer 3 (black squares). The dashed lines indicate matched pressure profiles using the Lanczos method with $\varepsilon_\text{sv} = 10^{-5}$. The solid lines show the upper and lower bounds for 20 inversions run using randomized 1-view with subspace re-use.}
\label{fig:9}    
\end{figure}

Figure \ref{fig:9} compares the history-matched pressure and enthalpy profiles found using the Lanczos method with those obtained using the 1-view method with subspace re-use. The 1-view method with subspace re-use gave the worst observation matches and was the only method to give an observation mismatch term outside the expected bounds \eqref{eq:obsmismespect}. The most noticeable difference between the methods is in the matched enthalpies for the high enthalpy Producer 1. Note that the oscillations that can be seen in the enthalpy matches for Producer 1 are probably due to numerical limitations of the model in handling flow of two-phase water in the reservoir. These oscillations make the inversion more challenging.

Inversion runs with the Lanczos method and the randomized methods not using subspace re-use gave good and very similar matches to the natural-state temperature observations. The methods applying subspace re-use tended to give slightly worse temperature matches though most of the matches were good.

\subsubsection{Increasing Truncation Faster}
\begin{figure}
\centering
\includegraphics[width=0.95\textwidth]{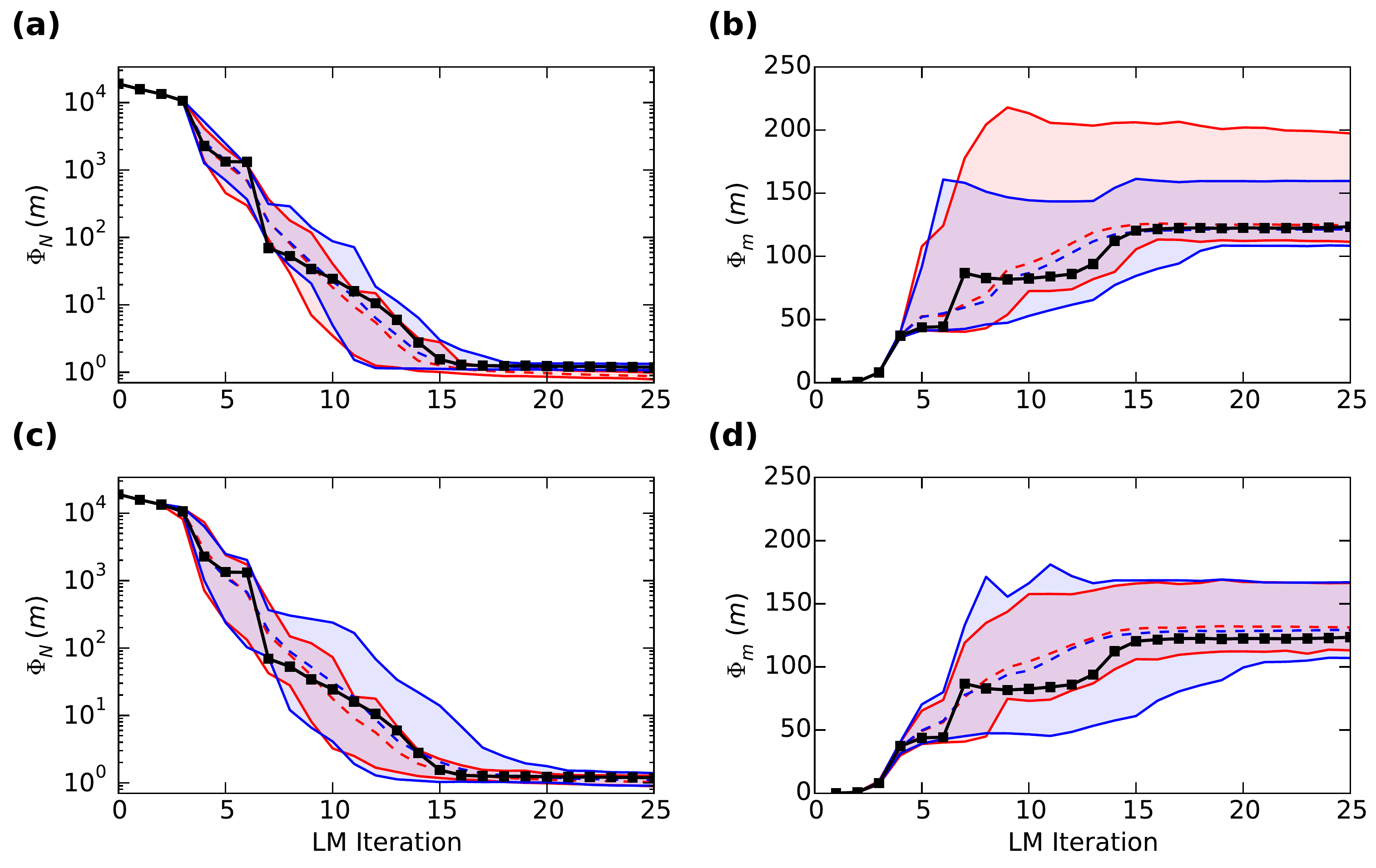}
\caption{Same as Fig. \ref{fig:7} but with the number of retained singular values increased linearly by five instead of two between LM iterations.}
\label{fig:10}       
\end{figure}

Additional inversions were run to look at the effects of increasing the number of singular values more rapidly when running the randomized TSVD-LM variants. Figure \ref{fig:10} compares the convergence behavior of the four randomized methods against the Lanczos method when increasing the number of retained singular values by five between iterations. The four randomized TSVD-LM methods were run 20 times and up to 25 LM iterations to generate the results shown in Fig. \ref{fig:10}. The results can be compared with those in Figs. \ref{fig:5}(c) and \ref{fig:5}(d).

The results in Fig. \ref{fig:10} again show that the convergence behavior of the inversions using the randomized methods is similar to inversions using the Lanczos method. However, there is greater variation in the model matches found using the randomized methods for this case compared to those shown in Fig. \ref{fig:7}. 

The increased variation is probably due to the rapid increase in the number of singular values, which increases the influence of the less accurately estimated singular triplets. This issue may be remedied by increasing the randomized oversampling and/or by considering the role of the number of retained singular values $p$, the LM damping factor $\gamma$ and the regularization weight $\mu$ in regularizing the inverse problem.

\section{Conclusions} \label{sec:conclusions}
The present study looked at applying randomized methods to improve the computational efficiency of inverting reservoir models using a modified LM approach based on the TSVD of a dimensionless sensitivity matrix. Randomized TSVD methods have not previously been applied to speeding up LM-based inversion of reservoir models. Previous methods have instead applied the iterative Lanczos method to find an approximate TSVD.

As discussed in this work the non-iterative nature of basic randomized TSVD methods allows for higher performance computing compared to classical iterative methods. Solution procedures which use randomized low-rank matrix approximation methods allow for simultaneous solution of the adjoint and direct problems used to form the TSVD of the dimensionless sensitivity matrix. 

Using a randomized 2-view method all the adjoint simulations can be run in parallel and likewise for the direct runs. Alternatively, the low-rank SVD approximation can be formed more efficiently by applying a random 1-view method, which enables solving all the adjoint problems concurrently with all the direct problems. This is in stark contrast to applying the standard Lanczos method which is inherently a serial process and therefore requires solving individual adjoint and direct problems one after the other. 

The proposed randomized methods were compared against the Lanczos approach for inverting permeability values of a simple, vertical slice, geothermal reservoir model. The randomized methods appear very promising as they resulted in an order of magnitude reduction in the computational time spent on running inversions compared with inversions using Lanczos iteration. Both the 1-view and 2-view methods converged to models with parameters and observation matches comparable with those found using the Lanczos approach. 

The 1-view method is especially appealing since a TSVD can in principle be evaluated in close to half the time needed by the 2-view method. However, the 2-view method may be preferable to the faster 1-view approach for certain problems, since the former can give better matches because of its superior accuracy. 

An alternative strategy could be to use the 1-view method for early LM iterations. A switch can then be made to a more accurate method, such as the 2-view one, at later stages of the inversion where the 1-view method can struggle.

The study also looked at improving the 1- and 2-view methods by re-using the subspace spanned by the singular vectors from previous LM iterations to initialize subsequent TSVD approximations. However, the randomized methods applying subspace re-use performed worse than the randomized methods without subspace re-use. The slightly simpler 1- and 2-view methods which do not apply subspace re-use are therefore preferable and merit further consideration for speeding up model inversion using the LM methodology.

Though the presented inversion approaches were demonstrated on a geothermal reservoir problem, the methodology can also be applied to related problems such as inversion of groundwater and petroleum reservoir models.

\appendix
\section{${\bm S}$ times a matrix ${\bm H}$} \label{sec:Direct}
Here we use TOUGH2 \citep{pruess2004,pruess1999tough2} as the reservoir simulator. Simulation of a geothermal system under production involves solving a natural state to find initial conditions for the production period. The discrete natural-state equations can be written as
\begin{linenomath*}
\begin{equation} \label{eq:natstate}
	{\bm f}_\text{st} ( {\bm u}_\text{st} , {\bm m} ) = {\bm 0} \, ,
\end{equation} \end{linenomath*}
where ${\bm u}_\text{st}$ are the primary variables for the natural state (e.g., pressures and temperatures). Due to the nonlinearity of the natural-state problem \eqref{eq:natstate} is solved using transient continuation. The production period that follows requires solving
\begin{linenomath*}
\begin{equation} \label{eq:tough2prod}
	{\bm f}^k ( {\bm u}^k , {\bm u}^{k-1} , {\bm m} , \Delta t^k ) = {\bm 0}
\end{equation} \end{linenomath*}
for every simulation time-step $\Delta t^k$ where $k$ is the time-step index and ${\bm u}$ holds the primary variables. The first production time-step is initialized with ${\bm u}^0 = {\bm u}_\text{st} ( {\bm m} )$.

The direct or forward method \citep{oliver2008,rodrigues2006adjoint} can be used to find ${\bm S}$ times a matrix ${\bm H} \in \mathbb{R}^{N_m \times s}$. ${\bm H}$ is a column vector ($s=1$) when running the Lanczos method. However, when applying randomized methods $s>1$ usually. For instance, $s=p+l$ when using the randomized 2-view method. 

Here we use the direct method based on \citet{oliver2008}, \citet{rodrigues2006adjoint} and the work presented in \citet{bjarkason2016} on applying adjoint and direct methods to geothermal simulations using TOUGH2. With the direct method we begin by solving for the natural state
\begin{linenomath*}
\begin{equation} \label{eq:dirnatstate}
	{\bm A}_\text{st} {\bm X}_\text{st} = - {\bm G}_\text{st} {\bm H}
\end{equation} \end{linenomath*}
to find ${\bm X}_\text{st}$ where
\begin{linenomath*}
\begin{equation}
	{\bm A}_\text{st} = \frac{\partial {\bm f}_\text{st}}{\partial {\bm u}_\text{st}} \quad \text{and} \quad {\bm G}_\text{st} = \frac{\partial {\bm f}_\text{st}}{\partial {\bm m}} \, .
\end{equation} \end{linenomath*}
Equation \eqref{eq:dirnatstate} is a linear problem with the number of right-hand sides equal to the number of columns in ${\bm H}$. For each iteration of the Lanczos method \eqref{eq:dirnatstate} has one right-hand side but \eqref{eq:dirnatstate} will generally be a multiple right-hand side problem when applying the randomized methods.

For the following production period we begin by solving for the first production time-step ($k = 1$)
\begin{linenomath*}
\begin{equation} \label{eq:dirprod1}
	{\bm A}^1 {\bm X}^1 = - {\bm G}^1 {\bm H} - {\bm B}^1 {\bm X}^0
\end{equation} \end{linenomath*}
to find ${\bm X}^1$ where
\begin{linenomath*}
\begin{equation} \label{eq:dirprodk}
	{\bm X}^0 = {\bm X}_\text{st} \, , \, 
    {\bm A}^k = \frac{\partial {\bm f}^k}{\partial {\bm u}^k} \, , \,
    {\bm B}^k = \frac{\partial {\bm f}^k}{\partial {\bm u}^{k-1}}
    \quad \text{and} \quad 
    {\bm G}^k = \frac{\partial {\bm f}^k}{\partial {\bm m}} \, .
\end{equation} \end{linenomath*}
For the time-steps that follow ($k > 1$) we solve
\begin{linenomath*}
\begin{equation}
	{\bm A}^k {\bm X}^k = - {\bm G}^k {\bm H} - {\bm B}^k {\bm X}^{k-1}
\end{equation} \end{linenomath*}
to obtain ${\bm X}^k$. Note that the dimensions of the linear problems in \eqref{eq:dirprod1} and \eqref{eq:dirprodk} are the same as \eqref{eq:dirnatstate}. With equations \eqref{eq:dirnatstate}, \eqref{eq:dirprod1} and \eqref{eq:dirprodk} for all $k > 1$ solved we find
\begin{linenomath*}
\begin{equation}
	{\bm S} {\bm H} = \frac{\mathrm{d} {\bm r}}{\mathrm{d} {\bm m}} {\bm H} = \frac{\partial {\bm r}}{\partial {\bm m}} {\bm H} + {\bm C}_\text{st} {\bm X}_\text{st} + \sum_{k=1}^{N_t^\text{pr}} {\bm C}^k {X^k} \, ,
\end{equation} \end{linenomath*}
where $N_t^\text{pr}$ is the number of production simulation time-steps and
\begin{linenomath*}
\begin{equation}
	{\bm C}_\text{st} = \frac{\partial {\bm r}}{\partial {\bm u}_\text{st}} \quad \text{and} \quad {\bm C}^k = \frac{\partial {\bm r}}{\partial {\bm u}^k} \, .
\end{equation} \end{linenomath*}

\section{${\bm S}^T$ times a matrix ${\bm H}$} \label{sec:Adjoint}
The adjoint method \citep{oliver2008,rodrigues2006adjoint} can be applied to estimate ${\bm S}^T$ times a matrix ${\bm H} \in \mathbb{R}^{N_d \times s}$. Like the direct method, the adjoint method requires solving linear problems at every simulation time-level. However, unlike the direct method the adjoint method tracks back in time.

Using the adjoint method (see \citet{bjarkason2016,oliver2008,rodrigues2006adjoint}) we begin by finding ${\bm Z}^{N_t^\text{pr}}$ for the final production time-step ($k=N_t^\text{pr}$)
\begin{linenomath*}
\begin{equation} \label{eq:adjlasttime}
	\left[ {\bm A}^{N_t^\text{pr}} \right]^T {\bm Z}^{N_t^\text{pr}} = - \left[ {\bm C}^{N_t^\text{pr}} \right]^T {\bm H} \, .
\end{equation} \end{linenomath*}

This is followed by tracking backwards in time solving for each preceding production simulation time
\begin{linenomath*}
\begin{equation} \label{eq:adjprod}
	\left[ {\bm A}^k \right]^T {\bm Z}^k = - \left[ {\bm C}^k \right]^T {\bm H}  - \left[ {\bm B}^{k+1} \right]^T {\bm Z}^{k+1}
\end{equation} \end{linenomath*}
to obtain ${\bm Z}^k$. Finally, we solve for the natural state
\begin{linenomath*}
\begin{equation} \label{eq:adjnatstate}
	\left[ {\bm A}_\text{st} \right]^T {\bm Z}_\text{st} = - \left[ {\bm C}_\text{st} \right]^T {\bm H}  - \left[ {\bm B}^{1} \right]^T {\bm Z}^{1}
\end{equation} \end{linenomath*}
to get ${\bm Z}_\text{st}$. The dimension of the linear problems in Eqs. (\ref{eq:adjlasttime}--\ref{eq:adjnatstate}) is the same as the linear problems solved for the direct method when applying the Lanczos method, the 2-view method or the 1-view method with $l_1=l_2$. The main difference is that the adjoint method requires tracking backward in time and the coefficient matrices of the linear problems (\ref{eq:adjlasttime}--\ref{eq:adjnatstate}) are transposes of the ones present in the direct method. 

After solving the above, the product of the transposed sensitivity matrix and the matrix ${\bm H}$ is found according to
\begin{linenomath*}
\begin{equation}
	{\bm S}^T {\bm H} = \left[ \frac{\mathrm{d} {\bm r}}{\mathrm{d} {\bm m}} \right]^T {\bm H} = \left[ \frac{\partial {\bm r}}{\partial {\bm m}} \right]^T {\bm H}  + \left[ {\bm G}_\text{st} \right]^T {\bm Z}_\text{st}  + \sum_{k=1}^{N_t^\text{pr}} \left[ {\bm G}^k \right]^T {\bm Z}^k  \, .
\end{equation} \end{linenomath*}

\section{Subspace Iteration Algorithm} \label{sec:appendixSI}
Algorithm \ref{alg:subiter} gives pseudocode for the subspace iteration method proposed by \citet{vogel1994} for finding a TSVD of a matrix ${\bm A}$. Their method requires evaluating ${\bm A}$ times a matrix and ${\bm A}^T$ times a matrix, just as the 1- and 2-view methods. The orthonormal initialization matrix ${\bm V }_0$ corresponds to the matrix ${\bm \Omega}$ in the 1- and 2-view methods.

\begin{algorithm}
  \caption{Subspace Iteration}\label{alg:subiter}
  \textbf{INPUT:} {Matrix ${\bm A} \in {\mathbb{R}}^{n_r \times n_c}$, integer $p > 0$, convergence tolerance $\varepsilon_\text{sv}$ and orthonormal matrix ${\bm V}^0 \in {\mathbb{R}}^{n_c \times p}$}.
  \begin{algorithmic}[1]
    \State Set $j=0$, {\ttfamily converged} = \textbf{False} and $\lambda^j_i = 0$, for $i=1, \, 2, \, \dots , \, p$.
    \State Evaluate ${\bm U}^j = {\bm A} {\bm V}^j$.  \label{algline:subiterU1}
    \While {(\textbf{not}  {\ttfamily converged})}
    	\State $j = j+1$.
        \State $\hat{{\bm C}} = {\bm A}^T {\bm U}^{j-1}$.
        \State Find the QR factors of $\hat{{\bm C}}$: $\hat{{\bm C}} = \hat{{\bm Q}} \hat{{\bm R}}$.
        \State Evaluate the SVD of $\hat{{\bm R}}$: $\hat{{\bm R}} = \hat{{\bm P}} \hat{{\bm \Lambda}} \hat{{\bm W}}^T$.
        \State ${\bm V}^j = \hat{{\bm Q}} \hat{{\bm P}}$.
        \State ${\bm U}^j = {\bm A} {\bm V}^j$.  \label{algline:subiterUj}
        \If {( Eq. \eqref{eq:lancConv} holds)}
        	\State {\ttfamily converged} = \textbf{True}.
        \EndIf
    \EndWhile
  \end{algorithmic}
\end{algorithm}

According to \citet{vogel1994} $\hat{\bm \Lambda}$ contains estimates of the $p$ largest singular values of the input matrix ${\bm A}$, ${\bm U}^j$ contains approximate left-singular vectors and ${\bm V}^j$ contains approximate right-singular vectors. However, $\hat{\bm \Lambda}$ actually contains estimates of the squares of the largest singular values on its diagonal, that is ${\bm \Lambda}_p \approx \hat{\bm \Lambda}^{1/2}$. Furthermore, the columns of ${\bm U}^j$ are not necessarily orthonormal as presented in Algorithm \ref{alg:subiter}, but this can be fixed with additional post-processing.

After a few iterations the left-singular vectors can be approximated by ${\bm U}_p \approx {\bm U}^j \hat{\bm \Lambda}^{-1/2}$. However, taking ${\bm U}_p$ as the left-singular vectors of ${\bm U}^j$, is a more accurate option.

Interestingly, a modern randomized 2-view method suggested by \citet{voronin2015} can be obtained by making minor changes to Algorithm \ref{alg:subiter}, assuming only one subspace iteration. The minor changes that need to be made to Algorithm \ref{alg:subiter} are to orthonormalize the matrix ${\bm U}^j$ in line \ref{algline:subiterU1}, replace ${\bm U}^j = {\bm A} {\bm V}^j$ on line \ref{algline:subiterUj} with ${\bm U}^j = \text{orth} ( {\bm U}^{j-1} ) \hat{\bm W}$ and introduce some oversampling. The resulting algorithm is given by Algorithm \ref{alg:2viewVoronin}. The difference between Algorithms \ref{alg:2view} and \ref{alg:2viewVoronin} is only in the way they use the matrix ${\bm B}^T$.

\begin{algorithm}
  \caption{Randomized Algorithm 2 in \citet{voronin2015}.}\label{alg:2viewVoronin}
  \textbf{INPUT:} {Matrix ${\bm A} \in {\mathbb{R}}^{n_r \times n_c}$ ($n_r \geq n_c$), integers $p > 0$ and $l \geq 0$.}\\
  \textbf{RETURNS:} {Approximate rank-$p$ SVD,  ${\bm U}_p {\bm \Lambda}_p {\bm V}_p^T$, of ${\bm A}$.}
  \begin{algorithmic}[1]		
    \State Generate a Gaussian random matrix ${\bm \Omega} \in {\mathbb{R}}^{n_c \times (p+l)}$.
    \State Form the matrix ${\bm Y} = {\bm A} {\bm \Omega}.$ \Comment{${\bm Y} \in {\mathbb{R}}^{n_r \times (p+l)}$}
    \State Find an orthonormal matrix ${\bm Q}\in {\mathbb{R}}^{n_r \times (p+l)}$, using QR factorization, such that ${\bm Y}={\bm Q} \tilde{\bm R}$.
    \State Evaluate the matrix ${\bm B}^T = {\bm A}^T {\bm Q}$. \Comment{${\bm B} \in {\mathbb{R}}^{(p+l) \times n_c}$}
    \State QR factorization, such that ${\bm B}^T = \hat{\bm Q} \hat{\bm R}$. \Comment{$\hat{\bm Q} \in {\mathbb{R}}^{n_c \times (p+l)}$, $\hat{\bm R} \in {\mathbb{R}}^{(p+l) \times (p+l)}$}  
    \State Calculate the SVD of the relatively small matrix $\hat{\bm R} = \hat{{\bm V}}_{p+l} {\bm \Lambda}_{p+l} \hat{\bm U}_{p+l}^T$ and truncate.
    \State Form the matrices ${\bm U}_p = {\bm Q} \hat{{\bm U}}_p$ and ${\bm V}_p = \hat{\bm Q} \hat{{\bm V}}_p$.
  \end{algorithmic}
\end{algorithm}

\acknowledgments
We thank N. Benjamin Erichson from the University of Washington, USA, for the introduction to the advantages of randomized SVD methods, which was a great inspiration for this study.
This work was supported by the NZ Ministry of Business, Innovation and Employment for funding this work through grant C05X1306. The first author was also kindly supported by Landsbankinn and an AUEA Braithwaite-Thompson Graduate Research Award.

\end{document}